\documentclass{amsart}

\usepackage[utf8]{inputenc}
\usepackage{pdflscape}
\usepackage{
 amsmath,
 amsxtra,
 amsthm,
 amssymb,
 etex,
 mathrsfs,
 }
\usepackage{hyperref}
\usepackage{xcolor}
\usepackage{longtable}
\usepackage{listings}
\usepackage{xcolor}
\usepackage{algpseudocode}
\usepackage{geometry}
\usepackage{algorithm}
\definecolor{codegreen}{rgb}{0,0.6,0}
\definecolor{codegray}{rgb}{0.5,0.5,0.5}
\definecolor{codepurple}{rgb}{0.58,0,0.82}
\definecolor{backcolour}{rgb}{0.95,0.95,0.92}

\lstdefinestyle{mystyle}{
    backgroundcolor=\color{backcolour},   
    commentstyle=\color{codegreen},
    keywordstyle=\color{magenta},
    numberstyle=\tiny\color{codegray},
    stringstyle=\color{codepurple},
    basicstyle=\ttfamily\footnotesize,
    breakatwhitespace=false,         
    breaklines=true,                 
    captionpos=b,                    
    keepspaces=true,                 
    numbers=left,                    
    numbersep=5pt,                  
    showspaces=false,                
    showstringspaces=false,
    showtabs=false,                  
    tabsize=2
}

\newtheorem{theorem}{Theorem}[section]
\newtheorem{lemma}[theorem]{Lemma}

\newtheorem{proposition}[theorem]{Proposition}
\newtheorem{corollary}[theorem]{Corollary}

\theoremstyle{definition}
\newtheorem{defn}[theorem]{Definition}
\newtheorem{remark}[theorem]{Remark}

\newcommand{\bd}{\begin{defn}}
\newcommand{\ed}{\end{defn}}
\newcommand{\bl}{\begin{lemma}}
\newcommand{\el}{\end{lemma}}
\newcommand{\bp}{\begin{proposition}}
\newcommand{\ep}{\end{proposition}}
\newcommand{\bt}{\begin{theorem}}
\newcommand{\et}{\end{theorem}}
\newcommand{\bc}{\begin{corollary}}
\newcommand{\ec}{\end{corollary}}
\newcommand{\br}{\begin{remark}}
\newcommand{\er}{\end{remark}}
\newcommand{\ba}{\begin{array}}
\newcommand{\ea}{\end{array}}
\newcommand{\bpf}{\begin{proof}}
\newcommand{\epf}{\end{proof}}

\newcommand{\Z}{\mathbb{Z}}
\newcommand{\Q}{\mathbb{Q}}

\newcommand{\Op}{\mathcal{O}}

\newcommand{\ze}{\zeta}

\DeclareMathOperator{\Gal}{Gal}
 \DeclareMathOperator{\rank}{rank}

\newcommand{\GL}{\mathrm{GL}}

\newcommand{\ot}{\otimes}
\newcommand{\ilim}{\displaystyle \mathop{\varinjlim}\limits}

\newcommand{\lra}{\longrightarrow}

  \DeclareFontFamily{U}{wncy}{}
  \DeclareFontShape{U}{wncy}{m}{n}{<->wncyr10}{}
  \DeclareSymbolFont{mcy}{U}{wncy}{m}{n}
  \DeclareMathSymbol{\sha}{\mathord}{mcy}{"58}

\begin{document}

\title{On even $K$-groups of rings of integers of real abelian fields}

\author{Meng Fai Lim}
\address{School of Mathematics and Statistics, Hubei Key Laboratory of Mathematical Sciences, Central China Normal University, Wuhan, 430079, P.R.China.}
\curraddr{}
\email{limmf@ccnu.edu.cn}
\thanks{M. F. Lim's research is partially supported by the Fundamental Research Funds for the Central Universities No. CCNU25JCPT031.
}

\author{Chao Qin}
\address{College of Mathematical Sciences,
Harbin Engineering University,
Harbin, 150001, P.R.China.}
\email{qinchao@hrbeu.edu.cn}
\thanks{Chao Qin's research is supported by the National Natural Science Foundation of China under Grant No. 12001546 and Heilongjiang Province under Grant No. 3236330122.}


\date{}

\dedicatory{}


\begin{abstract}
We present approaches for calculating the precise orders of the algebraic $K$-groups $K_{4n-2}(\mathcal{O}_K)$ for a totally real abelian $K$. Along the way, we also establish a formula connecting the order of $K_{4n-2}(\mathcal{O}_E)$ of a totally real $p$-elementary field $E$ to its intermediate cyclic $p$-degree fields. Additionally, we provide a compiled list of values for these $K$-groups.
\end{abstract}

\maketitle

\section{Introduction}

For a ring $R$, we let $K_{i}(R)$ denote the algebraic $K$-groups
of $R$ in the sense of Quillen \cite{Qui73b}. Thanks to the pioneering work of Quillen \cite{Qui73b}, Garland \cite{Gar}, and Borel \cite{Bo}, we now know that the even $K$-groups of the ring of integers of a number field are finite. However, these foundational results provide limited insight into the precise orders of these groups. It was the conjecture proposed by Birch and Tate that initially provided an approach to understanding the order of $K_2(\Op_F)$ by evaluating the Dedekind zeta function of the field $F$ at $s = -1$. This conjecture was later generalized by Lichtenbaum \cite{Lic72, Lic73} to include the higher even $K$-groups. Coates \cite{C72} then provided a crucial link by suggesting that this conjecture could be attacked using the main conjecture of Iwasawa theory \cite{Iw73}, an insight that has been pivotal in establishing the conjecture for totally real abelian fields and forms the backbone of our computational approaches. Since then, significant progress has been made towards establishing this conjecture, particularly notable in the case of a totally real abelian field.

Building on these conjectures, Browkin and his collaborators \cite{Brow00, Brow05, BrowGang} have gone a long way into computing the precise order of $K_2(\Op_F)$ for quadratic fields and specific classes of cubic fields (also see the work of Zhou \cite{Zh}). The primary objective of this paper is to continue this line of study by exploring the computation of the order of higher even $K$-groups, specifically $K_{4k-2}(\Op_F)$. More precisely, we discuss three different approaches for calculating the order of higher even $K$-groups, with each being particularly suited to certain types of number fields. The primary difference between these methods lies in the technique used to compute the special values of the Dedekind zeta function. Although most software tools offer built-in functions for numerical approximations of these values, our objective is to compute them exactly as fractions.

Our first approach applies in principle to all totally real abelian fields $F$. Leveraging the Artin formalism for $L$-functions, we establish a connection between the value $\zeta_F(1-2k)$ and a suitable product involving generalized Bernoulli numbers. Subsequently, we proceed to compute each of these generalized Bernoulli numbers individually. For this step, the associated Dirichlet character must be manually inputted from the LMFDB to ensure a precise fractional computation. We demonstrate this entire procedure in the body of the paper through a concrete example (see Appendix \ref{AppC}).
The second method is specifically for the real quadratic fields and is based on a formula of Siegel-Zagier. For this, we provide a simplified approach using Siegel-Zagier's formula and also offer a faster method to calculate $w_j(F)$, improving efficiency without compromising on obtaining exact values. This will be elaborated upon in Section \ref{Quadratic}.

The third approach is specifically tailored for a $p$-elementary totally real abelian field $E$, in which the Galois group $\Gal(E/\mathbb{Q})$ is isomorphic to $(\mathbb{Z}/p\mathbb{Z})^{\oplus n}$ for some prime $p$ and positive integer $n$. The key to this approach is the following theorem on which it relies.

\bt[Theorem \ref{mainthm}]
Let $p$ be a prime and $n$ an integer $\geq 2$. Suppose that $E$ is a totally real abelian extension of $\Q$ with Galois group $G=\Gal(E/\Q)\cong (\Z/p\Z)^{\oplus n}$. Let $K_1, K_2, \ldots, K_{\frac{p^n - 1}{p - 1}}$ denote all the $p$-degree extensions of $\mathbb{Q}$ contained in $E$. Then we have the following equality
\[ |K_{4k-2}(\Op_E)| = \frac{1}{|K_{4k-2}(\Z)|^{\frac{p^n-p}{p-1}}}\prod_{j=1}^{\frac{p^n-1}{p-1}}|K_{4k-2}(\Op_{K_j})|.\]
 \et

The proof of the theorem will be given in Section \ref{p-elementary section}. From a computational point of view, this result is rather advantageous to have. For instance, if we want to compute the size of $K_{4n-2}$-groups of the ring of integers of $\Q(\sqrt{2}, \sqrt{3}, \sqrt{5})$, it suffices to compute those for the intermediate quadratic fields $\Q(\sqrt{2})$, $\Q(\sqrt{3})$, $\Q(\sqrt{5})$, $\Q(\sqrt{6})$, $\Q(\sqrt{10})$, $\Q(\sqrt{15}), \Q(\sqrt{30})$, as well as for $\Q$. After obtaining these values, we can simply plug them into the formula provided by the aforementioned theorem.

\subsection*{Acknowledgement}

The collaboration began at a workshop at Harbin Institute of Technology in 2022. The authors extend their gratitude to Jun Wang and Yichao Zhang for organizing the workshop and for the insightful discussions that took place during the paper's preparation. Part of this research was conducted during the first author's multiple visits to Harbin Engineering University and the second author's visits to Central China Normal University. We would like to express our appreciation to both universities for their warm hospitality. We are particularly grateful to Daniel Delbourgo for his valuable insights, keen interest, and ongoing encouragement, all of which have been instrumental in advancing this work.

\section{$K$-groups and Dedekind zeta function} \label{review K}

We begin with a brief and quick review of the definition of the higher $K$-groups. Let $R$ be a ring with identity. For each integer $m \geq 1$, denote by $\GL_m(R)$ the group of invertible $m\times m$ matrices with entries in $R$. We then set $\GL(R) =\ilim_m\GL_m(R)$, where the transition map $\GL_m(R)\lra \GL_{m+1}(R)$ is given by
\[ A \mapsto \begin{pmatrix} A&0\\ 0&1 \\ \end{pmatrix}. \]
Let $\mathrm{BGL}(R)$ be the classifying space of the group $\GL(R)$, i.e., $\mathrm{BGL}(R)$ is an Eilenberg-MacLane space of type $(\mathrm{GL}(R),1)$ in the sense of \cite[Section 8.1]{Span}. Up to homotopy equivalence, this space is characterized by the property that it is path-connected with homotopy groups 
\[ \pi_n\big(\mathrm{BGL}(R)\big) \cong \left\{
                                  \begin{array}{ll}
                                    \GL(R), & \hbox{for $n=1$,} \\
                                    0, & \hbox{for $n\geq 2$.}
                                  \end{array}
                                \right.
\]
From the space $\mathrm{BGL}(R)$, one obtains a new space, denoted by $\mathrm{BGL}(R)^+$, via the $+$-construction of Quillen (see \cite{Qui73a, Qui73b}). The higher $K$-groups $K_n(R)$ are then defined by
\[ K_n(R):=\pi_n\big(\mathrm{BGL}(R)^+\big). \]
It's well-known that Quillen's construction recovers the classical $K_1$-groups of Bass and $K_2$-groups of Milnor (for instance, see \cite[Chapter IV]{WeiKbook}).

In this paper, we are interested in the $K$-groups of the ring $\Op_F$, where $\Op_F$ is the ring of integers of a number field $F$. As a start, we recall the following fundamental results of Quillen and Borel.

\bt \label{Borel}
The groups $K_{n}(\Op_F)$ are finitely generated for all $n\geq 1$. Furthermore, one has
 \[ \rank_{\Z} \big(K_{n}(\Op_F)\big) =\begin{cases} r_1(F)+r_2(F),  & \mbox{if $n \equiv 1$ (mod 4)}, \\
r_2(F), & \mbox{if $n \equiv 3$ (mod 4)},
\\
0, & \mbox{if $n$ is even}. \end{cases} \]
Here $r_1(F)$ (resp., $r_2(F)$) is the number of real embeddings (resp., number of pairs of complex embeddings) of the number field $F$.
\et

 \bpf
  Quillen \cite{Qui73b} was the first to establish that these $K$-groups are finitely generated. Subsequently, calculations of Borel \cite{Bo} confirmed the ranks of the $K$-groups as stated in the theorem. We should also mention that prior to the works of Qullen and Borel, the finiteness of $K_2(\Op_F)$ has been proven by Garland \cite{Gar}.
 \epf

Unfortunately, the results of Quillen and Borel do not provide a means to determine the exact order of the even $K$-groups $K_{2i}(\Op_F)$. It was only due to the remarkable insight of Birch, Tate and Lichtenbaum \cite{Lic72, Lic73} that one can hope to understand these orders via special values of Dedekind zeta function, a concept we will now outline briefly. Let $\ze_F(s)$ be the Dedekind zeta function of $F$. This function, $\ze_F(s)$, admits an analytic continuation to the whole complex plane, with the exception of a simple pole at $s = 1$. Consequently, it makes sense to speak of $\zeta_F(1-2k)$ for a positive integer $k$. Thanks to the collective gallant efforts of numerous mathematicians, we have the following.

\bt
Let $F$ be a totally real abelian number field of degree $r~ (=r_1(F))$. Then for every integer $k\geq 1$, we have
\[\zeta_F(1-2k)=(-1)^{kr}2^{r}\frac{|K_{4k-2}(\Op_F)|}{|K_{4k-1}(\Op_F)|}.\]
\et

\bpf
Lichtenbaum \cite{Lic72} first formulated this conjecture up to a power of $2$. Subsequently, the work of Coates \cite{C72} suggested that one might possibly attack this conjecture via the main conjecture of Iwasawa \cite{Iw73}. Building on this insight, Bayer and Neukirch \cite{BN} showed that the main conjecture of Iwasawa implies a cohomological version of Lichtenbaum's conjecture (for a detailed exposition of this cohomological version, readers are referred to \cite{BN}). Notably, this cohomological formulation is equivalent to the $K$-theoretical version, a connection established by the Quillen-Lichtenbaum conjecture. This conjecture is now a theorem, being a consequence of the groundbreaking work of Rost-Voevodsky (\cite{Vo}; see also Rognes-Weibel \cite{RW}). Prior to these developments, the main conjecture of Iwasawa has already been proven by by Mazur-Wiles \cite{MW} and Wiles \cite{Wiles}.
\epf

The value $|K_{4k-1}(\Op_F)|$ can be described rather easily. Let $\mu_\infty$ be the group of all the roots of unity of $\bar{F}$, where $\bar{F}$ is the algebraic closure of $F$. For an integer $j\geq 1$, we write $\mu_\infty^{\otimes j}$ for the $j$-fold tensor product of $\mu_\infty$ with $\Gal(\bar{F}/F)$ acting diagonally. Set $w_j(F)$ to be the order of $(\mu_\infty^{\otimes j})^{\Gal(\bar{F}/F)}$. The following gives a relation of $|K_{4k-1}(\Op_F)|$ in terms of these values.

\bt \label{K1 = w}
Let $F$ be a totally real abelian number field of degree $r$. Then for every integer $k\geq 1$, we have
\[ \big|K_{4k-1}(\Op_F)\big| =\begin{cases} 2^r w_{2k}(F),  & \mbox{if $k$ is odd}, \\
 w_{2k}(F), & \mbox{if $k$ is even}.
 \end{cases} \]
\et

\bpf
See \cite[Chap. VI, Theorem 9.5]{WeiKbook}.
\epf

Combining the above theorems, we have the following observation.

\bc \label{K2 w ze}
Let $F$ be a totally real abelian number field of degree $r$. Then for every integer $k\geq 1$, we have
\[ \big|K_{4k-2}(\Op_F)\big| =\begin{cases} (-1)^{r}w_{2k}(F)\zeta_F(1-2k),  & \mbox{if $k$ is odd}, \vspace{.1in}\\
 \displaystyle\frac{1}{2^r}w_{2k}(F)\zeta_F(1-2k), & \mbox{if $k$ is even}.
 \end{cases} \]
\ec

In principle, the values of $w_{2k}(F)$ can be determined rather easily (for instances, see \cite[Chap. VI, Propositions 2.2 and 2.3]{WeiKbook}).
Consequently, the main challenge lies in calculating the values of $\zeta_F(1-2k)$. For number fields of small degree, these can be computed via bulit-in functions in mathematical software programs. However, even in the case of a real quadratic field, the computed values are not exact when the discriminant of the real quadratic number field $F$ becomes large. In this context, computing $\zeta_F(-19)$ often leads to a loss of significant digits.  In the subsequent discussion of the paper, we will outline strategies to circumvent this issue. 
To begin with, we introduce the following approach that utilizes generalized Bernoulli numbers. 

Let $F$ be an abelian totally real number field with Galois group $G=\Gal(F/\Q)$. By Artin formalism, we have 
\[ \ze_F(1-2k)=\ze(1-2k)\prod_{\chi\neq \chi_0} L(\chi,1-2k)\]
where $\chi$ runs through all the nontrivial characters of $G$.
It's well-known that
\[ \ze(1-2k)=-\frac{B_{2k}}{2k} \quad\mathrm{and}\quad L(\chi, 1-2k) = -\frac{B_{2k,\chi}}{2k} \]
(for instance, see \cite[Theorem 4.2]{Wash}). 
The above therefore gives a way to compute the $\ze_F(1-2k)$ via generalized Bernoulli numbers. To see a specific example, we refer readers to the first listing in Appendix \ref{AppC}.

\section{Quadratic field} \label{Quadratic}
In this section, we shall describe a method (due to Siegel-Zagier) of computing the $L$-values for a totally real quadratic field.
To prepare for this, we need to introduce some further notations.
For a given integer $j\geq 0$ and an ideal $\mathfrak{a}$ of $\Op_F$, we define
\[ \sigma_j(\mathfrak{a}) = \sum_{\mathfrak{b}|\mathfrak{a}} |\Op/\mathfrak{b}|^j,\]
where the sum is taken over all nonzero ideals $\mathfrak{b}$ of $\Op_F$ that divide $\mathfrak{a}$.
In the special case where $\Op_F$ coincides with $\Z$, we shall simplify the notation to $\sigma(m) = \sigma(m\Z)$. Note that in this context, we have
 \[ \sigma_j(m) = \sum_{d|m} d^j,\]
where $d$ runs through all the positive divisors of $m$.
Furthermore, for integers $j,k\geq 1$, we set
\[ s_j^F(2k) = \sum_{ \substack{\nu\in\mathfrak{d}^{-1}\\ \nu\gg 0 \\ \mathrm{tr}(\nu)=j}}\sigma_{2k-1}\big((\nu)\mathfrak{d}\big).\]
Here $\mathfrak{d} = \mathfrak{d}_F$ is the different of $F$ and the sum is taken over all totally positive elements in $\mathfrak{d}$ with trace $j$. With these notation in hand, we can now state the following formula of Siegel.

\bt[Siegel] \label{Siegel formula}
Let $F$ be a totally real number field. Then for every integer $k\geq 1$, we have
\[ \ze_F(1-2k) = 2^{|F:\Q|}\sum_{j=1}^{r}b_j(2k|F:\Q|)s_j^F(2k), \]
where the numbers $b_j(2k|F:\Q|)$ are rational and depend only on $2k|F:\Q|$, and the integer $r$ is given by
\[ r =\begin{cases} \displaystyle\left[\frac{k|F:\Q|}{6}\right],  & \mbox{if $k|F:\Q|\equiv 1~(mod~ 6)$}, \vspace{.1in} \\
 \displaystyle\left[\frac{k|F:\Q|}{6}\right] +1, & \mbox{otherwise}.
 \end{cases} \]
\et

\br
Some of the values of Siegel coefficients $b_j(m)$ for $4\leq m\leq 40$ can be found in \cite[Table 1]{Za}. For the convenience of the readers, we provide an approach to computing these terms following \cite{Sie}. 

For $k=4,6,...$, we set
\[ G_k(z) = 1-\frac{2k}{B_k}\sum_{n=1}^\infty\sigma_{k-1}(n)q^n \]
and
\[\Delta = q\prod_{n=1}^\infty(1-q^n)^{24}.\]

We then define
\[ T_h = G_{12r-h+2}\Delta^{-r},\]
where \[r =
            \begin{cases}
              [\frac{h}{12}]+1, & \mbox{$h\not\equiv 2$ mod 12;} \\
              [\frac{h}{12}], & \mbox{$h\equiv 2$ mod 12.}
            \end{cases}
          \]

By \cite[P252, (11)]{Sie}, we have
\[ T_h = q^{-r} + c_{h,r-1}q^{-(r-1)}+\cdots + c_{h,1}q^{-1}+c_{h,0}+\cdots\]
for some $c_{h,i}\in \Q$. Furthermore, \cite[P253, Theorem 2]{Sie} tells us that $c_{h,0}\neq 0$. Hence it makes sense to write
\[b_{j}(h) := -c_{h,j}/c_{h,0}  \]
for $j=1,2,...,r$. These are the Siegel coefficients (see \cite[P254-255]{Sie} or \cite[P61, (26)]{Za}) that appears in Siegel's formula.

\er

For a real quadratic number field $K$, Zagier has expressed Siegel's formulas in terms of certain elementary functions which we now describe. For integers $j,m\geq 1$,  define
\[ e_j(m) = \sum_{\substack{b^2+4ac = m \\ a,c>0}}a^j. \]

We also denote by $\chi:=\chi_K$ the nontrivial character of $\Gal(K/\Q)$, and extend it to a function of $\Z$ in the usual way. Zagier's formula is then as follow (cf. \cite[(14),(16)]{Za}).

\bt[Zagier]
Let $K$ be a real quadratic field with discriminant $D$. Then for every integer $k\geq 1$, we have
\[ \ze_K(1-2k) = 4\sum_{j=1}^{{[k/3]+1}}b_j(4k) \sum_{m|j}\chi(m)m^{2k-1}e_{2k-1}\big((j/m)^2D\big). \]
\et

For a list of the formula for some values of $k$, we refer readers to \cite[Table 2]{Za}. A notable advantage of this formula lies in its ability to yield exact values of $\ze_K(1-2k)$ for real quadratic number fields with large discriminant. Moreover, by circumventing the need for $L$-value calculations in software like Pari or Magma, this formula potentially eliminates the dependency on GRH. 


We shall utilize the aforementioned formula to calculate the special values of the Dedekind zeta function of a real quadratic field. Furthermore, we introduce a faster method to compute $w_j(F)$, optimizing the process while maintaining the accuracy required to determine the sizes of the $K$ groups. This approach significantly improves computational efficiency in the quadratic case (see Listing 2 in Appendix \ref{AppC}).

\section{$p$-elementary extensions} \label{p-elementary section}

We are now in a position to prove the following theorem as mentioned in our introductory section.

\bt \label{mainthm}
Let $p$ be a fixed prime.
Let $E$ be a totally real Galois extension of $\Q$ with Galois group $G=\Gal(E/\Q)\cong (\Z/p\Z)^{\oplus n}$ with $n\geq 2$. Denote by $K_1, K_2,..., K_{\frac{p^n-1}{p-1}}$ all the $p$-degree extensions of $\Q$ contained in $E$. Then we have
\[ |K_{4k-2}(\Op_E)| = \frac{1}{|K_{4k-2}(\Z)|^{\frac{p^n-p}{p-1}}}\prod_{j=1}^{\frac{p^n-1}{p-1}}|K_{4k-2}(\Op_{K_j})|\]
for every positive integer $k$. \et

\bpf
For each $j=1,...,\frac{p^n-1}{p-1}$, we let $\chi_{j,r}$ be all the nontrivial character of $\Gal(K_j/\Q)$, where $r=1,...,p-1$. When viewed as characters of $\Gal(E/\Q)$, they are exactly all nontrivial characters of $\Gal(E/\Q)$.
By Artin formalism of $L$-functions, we have
\[ \ze_E(s) = \ze(s) \prod_{j=1}^{\frac{p^n-1}{p-1}} \prod_{r=1}^{p-1}L(\chi_{j,r},s). \]
On the other hand, we also have
\[ \ze_{K_j}(s) = \ze(s) \prod_{r=1}^{p-1}L(\chi_{j,r},s) \]
for each $j$. Hence we have
\[ \ze(s)^{\frac{p^n-p}{p-1}}\ze_E(s) = \prod_{j=1}^{\frac{p^n-1}{p-1}} \ze_{K_j}(s).\]
In view of the above equality and Corollary \ref{K2 w ze}, the proposition is reduced to proving the equality
\begin{equation}\label{w biquad full} w_{2k}(\Q)^{\frac{p^n-p}{p-1}} w_{2k}(E) =   \prod_{j=1}^{\frac{p^n-1}{p-1}} w_{2k}(K_j).
\end{equation}
 If $\ell$ is a prime, we write $w_j^{(\ell)}(F)$ for the order of $(\mu_{\ell^\infty}^{\otimes j})^{\Gal(\bar{F}/F)}$, where $\mu_{\ell^\infty}$ is the group of all the $\ell$-power roots of unity of $\bar{F}$. Plainly, one has $w_j(F) = \prod_\ell w_j^{(\ell)}(F)$. It therefore remains to show that
\begin{equation}\label{w biquad}
w_{2k}^{(\ell)}(\Q)^{\frac{p^n-p}{p-1}} w_{2k}^{(\ell)}(E) =   \prod_{j=1}^{\frac{p^n-1}{p-1}} w_{2k}^{(\ell)}(K_j)
\end{equation}
for every prime $\ell$.
We first consider the case when the prime $\ell$ is odd. Since $\Gal(E/\Q)$ is not cyclic, we have either $E\cap \Q(\mu_\ell)=\Q$ or $E\cap \Q(\mu_\ell) = K_i$ for some unique $i$.

Suppose that $E\cap \Q(\mu_\ell)=\Q$. Then we have $|L(\mu_\ell):L|= \ell-1$ for $L= E, K_j, \Q$. If $2k$ is not divisible by $\ell-1$, it follows from \cite[Chap.\ VI, Proposition 2.2(c)]{WeiKbook} that $w_{2k}^{(\ell)}(L) =1$ for $L= E, K_j, \Q$, and so equality (\ref{w biquad}) is immediate in this case. If $2k$ is divisible by $\ell-1$, then \cite[Chap.\ VI, Proposition 2.2(c)]{WeiKbook} tells us that $w_{2k}^{(\ell)}(L) = \ell^{1+b}$ for $L= E, K_j, \Q$, where $b$ is the highest power of $\ell$ dividing $2k$. This again verifies the equality (\ref{w biquad}).

Now, without loss of generality, suppose that $E\cap \Q(\mu_\ell) = K_1$. In other words, $K_1$ is contained in $\Q(\mu_\ell)$ with $\ell\equiv 1$ (mod $2p$). If $2k$ is not divisible by $(\ell-1)/p$,  it follows from \cite[Chap.\ VI, Proposition 2.2(c)]{WeiKbook} that $w_{2k}^{(\ell)}(L) =1$ for $L= E, K_j, \Q$, thus verifying the equality (\ref{w biquad}). In the event that $2k$ is divisible by $\ell-1$, a similar argument as in the previous paragraph yields the required equality (\ref{w biquad}). Therefore, it remains to consider the case where $2k$ is divisible by $(\ell-1)/p$ but not divisible by $\ell-1$. In this case, one can directly verify that $w_{2k}^{(\ell)}(E) =  w_{2k}^{(\ell)}(K_1) = \ell^{1+b}$ and $w_{2k}^{(\ell)}(\Q) =  w_{2k}^{(\ell)}(K_j)=1$ for $j\geq 2$. From this, we see that the equality (\ref{w biquad}) is satisfied.

We now come to the situation when $\ell=2$. We first consider the case $\sqrt{2}\notin E$. Under this said assumption, the field $L(\sqrt{-1})$ does not contain any primitive $8$th root of unity. Thus, it follows from \cite[Chap.\ VI, Proposition 2.3(c)]{WeiKbook} that $w_{2k}^{(2)}(L) = 2^{2+d}$ for $L= E, K_j, \Q$, where $d$ is the highest power of $2$ dividing $2k$. Plainly, the equality (\ref{w biquad}) is satisfied in this case. Now suppose that $\sqrt{2}\in E$. In particular, we must then have $p=2$. Upon relabeling, we may assume $K_1= \Q(\sqrt{2})$. Note that the fields $K_1(\sqrt{-1})$ and $E(\sqrt{-1})$ will now contain a primitive $8$th root of unity but the fields $\Q(\sqrt{-1})$ and $K_j(\sqrt{-1})$ (for $j\geq 2$) do not. Hence, from \cite[Chap.\ VI, Proposition 2.2(c,d)]{WeiKbook} it follows that we have $w_{2k}^{(2)}(E) =  w_{2k}^{(2)}(K_1) = 2^{3+d}$ and $w_{2k}^{(2)}(\Q) =  w_{2k}^{(2)}(K_j)=2^{2+d}$ for $j\geq 2$. Plugging these values into the equality (\ref{w biquad}), we see that the said equality holds.

The proof of the theorem is therefore complete.
 \epf

  Note that the asserted equality in the preceding proposition is not true if one remove the ``totally real" hypothesis. Indeed, for an imaginary biquadratic field, Guo and Qin has shown that there might be an extra factor of a power of 2 (see \cite[Theorem 3.5, Example 3.10]{GuoQin}).

From a computational point of view, Theorem \ref{mainthm} is rather advantageous to have, as it streamlines the process of computing the $K$-group for an elementary $p$-extension of $\mathbb{Q}$ by reducing it to the computation of the $K$-group for a cyclic $p$-degree extension of $\mathbb{Q}$. Regarding the latter task, and in light of Corollary \ref{K2 w ze}, our focus shifts to determining the values of $w_{2k}$ and the relevant special values. For this, we present the following useful proposition.

\bp\label{ww}
Let $p$ be a prime and $k\geq 1$. Suppose that $K$ is a totally real number field which is a cyclic extension of $\Q$ of degree $p$. Then the following statements are valid.
\begin{enumerate}
  \item[$(i)$] $w_{2k}(\Q(\sqrt{2}))=  \displaystyle 2^{4+v_2(k)}\prod_{\substack{l~\mbox{\footnotesize odd prime}\\ l-1\mid 2k}}l^{1+v_l(k)}=  \displaystyle \prod_{\substack{l~\mbox{\footnotesize  prime}\\ l-1\mid 2k}}l^{1+v_l(8k)}$.
  \item[$(ii)$] If the prime $p$ is odd and $K\subseteq \Q(\ze_{p^2})$, then
 \[ w_{2k}(K)=2^{3+v_2(k)}p^{(2+v_p(k))\delta}\prod_{\substack{l\neq 2,p \\ l-1\mid 2k}}l^{1+v_l(k)}  = \prod_{\substack{l~\mbox{\footnotesize  prime} \\ l-1\mid 2k}}l^{1+v_l(4pk)},\]
where $\delta =\left\{
                                       \begin{array}{ll}
                                         1, & \hbox{if $p-1\mid 2k;$} \\
                                         0, & \hbox{otherwise.}
                                       \end{array}
                                     \right.$
  \item[$(iii)$] If $K$ is contained in $\Q(\ze_q)$ for some odd prime $q$, then
\[ w_{2k}(K)=2^{3+v_2(k)}q^{(1+v_q(k))\varepsilon}\prod_{\substack{l\neq 2,q \\ l-1\mid 2k}}l^{1+v_l(k)},\]
where $\varepsilon =\left\{
                                       \begin{array}{ll}
                                         1, & \hbox{if $\frac{q-1}{p}\mid 2k;$} \\
                                         0, & \hbox{otherwise.}
                                       \end{array}
                                     \right.$
\item[$(iv)$] For other $K$'s not covered in $(i)-(iii)$, one always has
\[ w_{2k}(K)=2^{3+v_2(k)}\prod_{\substack{l\neq 2 \\ l-1\mid 2k}}l^{1+v_l(k)}= \prod_{\substack{l~\mbox{\footnotesize  prime} \\ l-1\mid 2k}}l^{1+v_l(4k)}.\]
\end{enumerate}
\ep

\bpf
Indeed, for an odd prime $p$, it follows from \cite[Chap. VI, Proposition 2.2]{WeiKbook} that
\[ \mu_{p^{\infty}}^{\ot i}(K) = \left\{
                                   \begin{array}{ll}
                                     \mu_{p^{a(K)+v_p(i)}}^{\ot i}, & \hbox{if $i \equiv 0$ mod $|K(\mu_p):K|$,} \\
                                     1, & \hbox{if $i \not\equiv 0$ mod $|K(\mu_p):K|$,}
                                   \end{array}
                                 \right.
 \]
 where
 $a(K)$ is the largest integer such that $K(\mu_p)$ contains a primitive $p^{a(K)}$th root of unity. For $p=2$ and even $i$, an application of \cite[Chap. VI, Proposition 2.3(c)]{WeiKbook} tells us that
\[ \mu_{2^{\infty}}^{\ot i}(K) =                                      \mu_{2^{c(K)+v_2(i)}}^{\ot i}
 \]
 where
 $c(K)$ is the largest integer such that $K(\sqrt{-1})$ contains a primitive $2^{c(K)}$th root of unity. (Note that our number field $K$ is totally real and so is an exceptional one in the sense of the proposition loc. cit.) The conclusions of the proposition now follow from the above two observations and a case-by-case analysis. 
\epf

For the code designed for the computation of even $K$-groups within the context of a $p$-elementary extension, see Listing 3 in Appendix \ref{AppC}.

\section{Periodicity of $p$-rank} \label{periodicity section}
Let $K$ be an real quadratic field. Recall that for integers $j,m\geq 1$, we have defined
\[ e_j(m) = \sum_{\substack{b^2+4ac = m \\ a,c>0}}a^j. \] 
Plainly, one has $e_1(D)\equiv e_3(D)$ (mod 3). On the other hand, we have \[\big|K_2(\Op_K)\big| = \left\{
            \begin{array}{ll}
              \displaystyle\frac{4}{5} e_1(8), & \hbox{$K=\Q(\sqrt{2})$,} \vspace{.1in} \\
             2e_1(5), & \hbox{$K=\Q(\sqrt{5})$,} \vspace{.1in}\\
              \displaystyle\frac{2}{5} e_1(D), & \mbox{otherwise.}
            \end{array} \right. \quad \quad\big|K_6(\Op_K)\big| = \left\{
            \begin{array}{ll}
              e_3(8), & \hbox{$K=\Q(\sqrt{2})$,} \vspace{.1in} \\
              \displaystyle\frac{1}{2} e_3(D), & \mbox{otherwise.}
            \end{array}
          \right. \] 
Therefore, it follows that $\big|K_2(\Op_K)\big|$ is divisible by $3$ if and only if $\big|K_6(\Op_K)\big|$ is divisible by $3$. Indeed, we shall see that there is a general reasoning underlying this observation which will be elucidated in the subsequent theorem.

\bt \label{periodicity of p-rank}
Let $F$ be a number field. Then we have
\[ \rank_{\Z/p\Z}\big( K_{2i}(\Op_F)\big) = \rank_{\Z/p\Z}\big( K_{2i'}(\Op_F)\big), \]
whenever $i \equiv i'$ $($mod $|F(\mu_p):F|)$.
\et

\bpf
 Although this fact might be well-known among experts but for the convenience of the readers, we shall provide a brief outline of the proof here. By the work of Rost and Voevodsky \cite{Vo}, there is an identification
\[ K_{2k}(\Op_F)/p \cong H^2\big(\Gal(F_{S_p}/F)), \mu_p^{\otimes(k+1)}\big),\]
where $F_{S_p}$ is the maximal algebraic extension of $F$ unramified outside the set of primes of $F$ above $p$.

If $j \equiv 0$ mod $[F(\mu_p) : F]$, then the Galois group $\Gal(F_{S_p}/F)$ acts trivially on $\mu_p^{\otimes j}$. Therefore, it follows that 
\[
  \begin{split}
    K_{2k}(\Op_F)/p  & \cong H^2\big(\Gal(F_{S_p}/F)), \mu_p^{\otimes(k+1)}\big) \\
       &\cong H^2\big(\Gal(F_{S_p}/F)), \mu_p^{\otimes(k'+1)}\big)\otimes \mu_p^{\otimes(k-k')}\\
       &\cong K_{2k'}(\Op_F)/p\otimes \mu_p^{\otimes(k-k')},
   \end{split}\]
  whenever $k \equiv k'$ mod $[F(\zeta_p) : F]$.
  Consequently, the groups $K_{2k}(\Op_F)/p$ and $K_{2k'}(\Op_F)/p$ have the same rank over $\Z/p\Z$.
\epf

 We return to the context of a real quadratic field $K$. Proposition \ref{periodicity of p-rank} then tells us that the $3$-rank
\[ r_3\big( K_{4k-2}(\Op_K)\big)\]
is a constant function in term of $k$. A consequence of this is the following.

\bc
Let $K$ be a real quadratic field with discriminant $D$. Denote by $\chi$ the notrivial character of $\Gal(K/\Q)$. Then the following statements are equivalent.\begin{enumerate}
             \item[$(1)$] $e_1(D)$ is divisible by $3$.
             \item[$(2)$]  $e_3(D)$ is divisible by $3$.
             \item[$(3)$]   $e_5(4D) + \big(5\chi(2) + 6\big)e_5(D)$ is divisible by $9$.
             \item[$(4)$]  $e_7(4D) + 19\chi(2)e_7(D)$ is divisible by $27$.
             \item[$(5)$]   $e_9(4D) + \big(8\chi(2)+3\big)e_9(D)$ is divisible by $9$.
 \item[$(6)$] $e_{11}(9D)$ is divisible by $3$.
\item[$(7)$] $e_{13}(9D)$ is divisible by $3$.
\item[$(8)$] $e_{15}(9D)$ is divisible by $3$.
           \end{enumerate}
\ec

\br
We remark that the list in the preceding corollary is far from exhaustive and goes on. One can of course apply the above discussion for other primes. For instances, for the case $p=5$, then the following divisibility statements are equivalent.
\begin{enumerate}
             \item[$(1)$] $e_1(D)$ is divisible by $25$.
             \item[$(2)$]   $e_5(4D) + \big(7\chi(2) -1\big)e_5(D)$ is divisible by $25$.
            \item[$(3)$]   $e_9(4D) + \big(8\chi(2)+3\big)e_9(D)$ is divisible by $25$.
           \item[$(4)$] $e_{13}(9D)$ is divisible by $5$.
           \end{enumerate}

One might naturally wonder whether these divisibility implications can be directly explained through the lens of these power sums, although we will not explore this particular subject in the current paper.
\er

We end with another possible application of Theorem \ref{periodicity of p-rank}. It is a natural question to ask whether $K_{2i}(\Op_F)[p^\infty]$ is cyclic for a given prime $p$. The following corollary gives a sufficient condition for verifying this.

\bc
Let $F$ be a number field. Suppose that $\big| K_{2i_0}(\Op_F)[p^\infty]\big| = p$ for some $i_0$. Then we have
\[ r_p\big( K_{2i}(\Op_F)\big) =1, \]
whenever $i \equiv i_0$ $($mod $|F(\mu_p):F|)$. In other words, $K_{2i}(\Op_F)[p^\infty]$ is cyclic for $i \equiv i_0$ $($mod $|F(\mu_p):F|)$.
\ec

\newpage
\appendix
\section{Algorithms} \label{AppA}
This appendix provides algorithms for computing the size of higher even $K$-groups of ring of integers of totally real number fields. The algorithms are implemented in SageMath and are divided into two categories: those for quadratic fields, including Siegel formula-related computations and multi-quadratic field cases, and those for general number fields with Galois group isomorphic to \((\mathbb{Z}/p\mathbb{Z})^n\). 
\subsection{Quadratic Field Case}
This section presents algorithms tailored for quadratic fields, including foundational computations, Siegel formula-related functions, and a method for combining $K$-group sizes across multiple quadratic fields. These algorithms compute the size of higher even $K$-groups and related quantities using modular forms, zeta functions, and prime products.

\begin{algorithm}[H]
\caption{Computes the size of $K_{n}(\mathbb{Z})$}
\begin{algorithmic}[1]
\State \textbf{Input}: Integer \( n \)
\State \textbf{Output}: Absolute value of constant
\Function{KZ}{$n$}
    \If{\( n \mod 8 = 2 \)}
        \State \Comment{Case: \( n \equiv 2 \pmod{8} \)}
        \State Compute \( k \gets (n - 2)/8 \)
        \State Compute \( s \gets 2k + 1 \)
        \State Compute \( c \gets \text{bernoulli}(2s) \cdot w(\mathbb{Q}, 1, 2s) / (4s) \)
        \State \Return \( |2c| \)
    \ElsIf{\( n \mod 8 = 6 \)}
        \State \Comment{Case: \( n \equiv 6 \pmod{8} \)}
        \State Compute \( k \gets (n - 6)/8 \)
        \State Compute \( s \gets 2k + 2 \)
        \State Compute \( c \gets \text{bernoulli}(2s) \cdot w(\mathbb{Q}, 1, 2s) / (4s) \)
        \State \Return \( |c| \)
    \EndIf
\EndFunction
\end{algorithmic}
\end{algorithm}

\vspace{0.5cm}

\begin{algorithm}[H]
\caption{Compute $w_j(F)$}
\begin{algorithmic}[1]
\State \textbf{Input}: Quadratic field \( K \), integer \( i \)
\State \textbf{Output}: Product of primes
\Function{w}{$K, i$}
    \State Initialize empty list \( \text{list1} \)
    \State \Comment{Collect primes where \( i/(\ell-1) \in \mathbb{Z} \)}
    \For{\( \ell \) in \( \text{primes}(2, i+2) \)}
        \If{\( i/(\ell-1) \in \mathbb{Z} \)}
            \State Append \( \ell \) to \( \text{list1} \)
        \EndIf
    \EndFor
    \State Filter \( \text{list1} \gets \{ \ell \in \text{list1} \mid \ell \neq 2 \} \)
    \State Compute \( k \gets 2i \)
    \State Initialize empty list \( \text{list2} \)
    \State \Comment{Collect additional primes}
    \For{\( \ell \) in \( \text{primes}(2, 2i+2) \)}
        \If{\( k/(\ell-1) \in \mathbb{Z} \) and \( i/(\ell-1) \not\in \mathbb{Z} \)}
            \State Append \( \ell \) to \( \text{list2} \)
        \EndIf
    \EndFor
    \State Initialize \( \text{prod} \gets 1 \)
    \State \Comment{Compute product for list1 primes}
    \For{\( \ell \in \text{list1} \)}
        \If{\( i/\ell \not\in \mathbb{Z} \)}
            \State \( \text{prod} \gets \text{prod} \cdot \ell \)
        \Else
            \State \( \text{prod} \gets \text{prod} \cdot \ell^{\text{valuation}(i, \ell) + 1} \)
        \EndIf
    \EndFor
    \State \( \text{prod} \gets \text{prod} \cdot 2^{\text{valuation}(i, 2) + 2} \)
    \State \Comment{Adjust for specific quadratic fields}
    \For{\( \ell \in \text{list2} \)}
        \If{\( K = \mathbb{Q}(\sqrt{\ell}) \)}
            \State \( \text{prod} \gets \text{prod} \cdot \ell^{\text{valuation}(i, \ell) + 1} \)
        \EndIf
    \EndFor
    \If{\( K = \mathbb{Q}(\sqrt{2}) \) or \( K = \mathbb{Q}(\sqrt{8}) \)}
        \State \( \text{prod} \gets 2 \cdot \text{prod} \)
    \EndIf
    \State \Return \( \text{prod} \)
\EndFunction
\end{algorithmic}
\end{algorithm}

\vspace{0.5cm}

\begin{algorithm}[H]
\caption{Compute $T_h=G_{12r-h+2}\Delta^{-r}$}
\begin{algorithmic}[1]
\State \textbf{Input}: Integer \( h \)
\State \textbf{Output}: Modular form \( T \)
\Function{T}{$h$}
    \If{\( h \mod 12 = 2 \)}
        \State Compute \( r \gets \lfloor h/12 \rfloor \)
    \Else
        \State Compute \( r \gets \lfloor h/12 \rfloor + 1 \)
    \EndIf
    \State Compute \( k \gets 12r - h + 2 \)
    \If{\( k = 0 \)}
        \State Compute \( T \gets \Delta^{-r} \)
    \Else
        \State Compute Eisenstein series \( G \gets \text{eisenstein\_series\_qexp}(k, r+1, \text{normalization='constant'}) \)
        \State Compute \( T \gets G \cdot \Delta^{-r} \)
    \EndIf
    \State \Return \( T + O(q) \)
\EndFunction
\end{algorithmic}
\end{algorithm}

\vspace{0.5cm}

\begin{algorithm}[H]
\caption{Compute Siegel coefficients $b_j(h)$}
\begin{algorithmic}[1]
\State \textbf{Input}: Integers \( l \), \( h \)
\State \textbf{Output}: Ratio of modular form coefficients
\Function{b}{$l, h$}
    \State Compute modular form \( T \gets T(h) \)
    \State Extract coefficients \( C \gets T.\text{coefficients}() \)
    \State Compute length \( L \gets |C| \)
    \State \Return \( -C[L - l - 1] / C[L - 1] \)
\EndFunction
\end{algorithmic}
\end{algorithm}

\vspace{0.5cm}

\begin{algorithm}[H]
\caption{Compute Siegel Sum}
\begin{algorithmic}[1]
\State \textbf{Input}: Integer \( D \), exponent \( j \)
\State \textbf{Output}: Sum of divisors raised to power \( j \)
\Function{eSiegel}{$D, j$}
    \State Initialize \( \text{haha} \gets 0 \)
    \State Initialize \( \text{bla} \gets 0 \)
    \If{\( D/4 \in \mathbb{Z} \)}
        \State \Comment{Sum over divisors of \( D/4 \)}
        \For{each \( i \in \text{divisors}(\lfloor D/4 \rfloor) \)}
            \State \( \text{haha} \gets \text{haha} + i^j \)
        \EndFor
    \EndIf
    \State \Comment{Sum over divisors for quadratic residues}
    \For{\( b \) in \( [1, \lfloor \sqrt{D} \rfloor] \)}
        \If{\( (D - b^2)/4 \in \mathbb{Z} \)}
            \For{each \( i \in \text{divisors}(\lfloor (D - b^2)/4 \rfloor) \)}
                \State \( \text{bla} \gets \text{bla} + i^j \)
            \EndFor
        \EndIf
    \EndFor
    \State \Return \( \text{haha} + 2 \cdot \text{bla} \)
\EndFunction
\end{algorithmic}
\end{algorithm}

\vspace{0.5cm}

\begin{algorithm}[H]
\caption{Compute Weighted Sum with Kronecker Symbol}
\begin{algorithmic}[1]
\State \textbf{Input}: Quadratic field \( F \), integers \( l \), \( m \)
\State \textbf{Output}: Weighted sum
\Function{S}{$F, l, m$}
    \State Compute discriminant \( D \gets F.\text{discriminant}() \)
    \State Initialize \( \text{sum} \gets 0 \)
    \State \Comment{Sum over divisors of \( l \)}
    \For{each \( j \in \text{divisors}(l) \)}
        \State Compute \( \text{sum} \gets \text{sum} + \text{kronecker\_symbol}(D, j) \cdot j^{2m-1} \cdot \text{eSiegel}((l/j)^2 \cdot D, 2m-1) \)
    \EndFor
    \State \Return \( \text{sum} \)
\EndFunction
\end{algorithmic}
\end{algorithm}

\vspace{0.5cm}

\begin{algorithm}[H]
\caption{Compute Zeta Function Value for Quadratic Field}
\begin{algorithmic}[1]
\State \textbf{Input}: Discriminant \( D \), complex number \( s \)
\State \textbf{Output}: Zeta function value
\Function{Zeta}{$D, s$}
    \State Create quadratic field \( F \gets \mathbb{Q}(\sqrt{D}) \)
    \State Set \( r \gets 2 \)
    \State Compute \( k \gets (1 - s)/2 \)
    \If{\( kr \mod 6 = 1 \)}
        \State Compute \( c \gets \lfloor kr/6 \rfloor \)
    \Else
        \State Compute \( c \gets \lfloor kr/6 \rfloor + 1 \)
    \EndIf
    \State Initialize \( \text{sum} \gets 0 \)
    \State \Comment{Sum over coefficients}
    \For{\( j \) in \( [1, c] \)}
        \State Compute \( \text{sum} \gets \text{sum} + b(j, 2kr) \cdot S(F, j, k) \)
    \EndFor
    \State \Return \( 2^r \cdot \text{sum} \)
\EndFunction
\end{algorithmic}
\end{algorithm}

\vspace{0.5cm}

\begin{algorithm}[H]
\caption{Compute size of $K$-group of Quadratic Field}
\begin{algorithmic}[1]
\State \textbf{Input}: Quadratic field \( K \), integer \( n \)
\State \textbf{Output}: Rounded K-group size
\Function{KSize}{$K, n$}
    \State Compute degree \( r \gets K.\text{degree}() \)
    \State Compute discriminant \( D \gets K.\text{discriminant}() \)
    \If{\( n \mod 2 = 0 \)}
        \State Compute \( k \gets (n + 2)/4 \)
        \If{\( k \mod 2 = 0 \)}
            \State \Comment{Even \( k \)}
            \State Compute \( SS \gets (1/2^r) \cdot w(K, 2k) \cdot \text{Zeta}(D, 1-2k) \)
        \Else
            \State \Comment{Odd \( k \)}
            \State Compute \( SS \gets (-1)^r \cdot w(K, 2k) \cdot \text{Zeta}(D, 1-2k) \)
        \EndIf
    \Else
        \State Compute \( k \gets (n + 1)/4 \)
        \If{\( k \mod 2 = 0 \)}
            \State \Comment{Even \( k \)}
            \State Compute \( SS \gets w(K, 2k) \)
        \Else
            \State \Comment{Odd \( k \)}
            \State Compute \( SS \gets 2^r \cdot w(K, 2k) \)
        \EndIf
    \EndIf
    \State \Return \( \text{round}(SS) \)
\EndFunction
\end{algorithmic}
\end{algorithm}

\vspace{0.5cm}

\begin{algorithm}[H]
\caption{Compute the size of $K$-group Size for multi-quadratic Fields}
\begin{algorithmic}[1]
\State \textbf{Input}: List of quadratic fields \( \text{list1} \), list of conductors \( \text{list2} \), integer \( n \)
\State \textbf{Output}: Combined K-group size
\Function{MultiSize}{$\text{list1}, \text{list2}, n$}
    \State Initialize \( \text{prod} \gets 1 \)
    \State Compute degree \( p \gets \text{list1}[0].\text{degree}() \)
    \State Set \( k \gets 3 \) \Comment{Assumes Galois group degree}
    \State \Comment{Compute product of K-group sizes}
    \For{\( i \) in \( [0, |\text{list1}|-1] \)}
        \State Compute \( \text{prod} \gets \text{prod} \cdot \text{KSize}(\text{list1}[i], \text{list2}[i], n) \)
    \EndFor
    \State Compute denominator \( \text{denominator} \gets \text{KZ}(n)^{(p^k - p)/(p - 1)} \)
    \State \Return \( \text{prod} / \text{denominator} \)
\EndFunction
\end{algorithmic}
\end{algorithm}
\textbf{Purpose}: Computes the combined K-group size for a list of quadratic fields.

\subsection{General Number Field Case}
This section presents algorithms for computing K-group sizes for general number fields, particularly those with Galois group \((\mathbb{Z}/p\mathbb{Z})^{\oplus n}\). These algorithms handle subfield checks, prime products, zeta functions, and K-group computations for arbitrary degrees.

\begin{algorithm}[H]
\caption{Check if \( K \) is a Subfield of \( L \)}
\begin{algorithmic}[1]
\State \textbf{Input}: Number fields \( K \), \( L \)
\State \textbf{Output}: Boolean indicating if \( K \subseteq L \)
\Function{is\_subfield}{$K, L$}
    \State Compute embeddings \( \text{Emb}(K, L) \) of \( K \) into \( L \)
    \State \Return \( |\text{Emb}(K, L)| > 0 \)
\EndFunction
\end{algorithmic}
\end{algorithm}

\vspace{0.5cm}

\begin{algorithm}[H]
\caption{Compute $w_j(F)$ for a random totally real field $F$, based on Proposition \ref{ww}}
\begin{algorithmic}[1]
\State \textbf{Input}: Number field \( K \), conductor \( \text{cond} \), integer \( i \)
\State \textbf{Output}: Product of prime powers
\Function{w}{$K, \text{cond}, i$}
    \State Compute degree \( r \gets \deg(K) \)
    \State Initialize empty list \( \text{list1} \)
    \State \Comment{Collect primes where \( i \mod (\ell - 1) = 0 \)}
    \For{\( \ell \) in \( \text{prime\_range}(2, i+2) \)}
        \If{\( i \mod (\ell - 1) = 0 \)}
            \State Append \( \ell \) to \( \text{list1} \)
        \EndIf
    \EndFor
    \State Initialize \( \text{prod} \gets 1 \)
    \If{\( K = \mathbb{Q}(\sqrt{2}) \)}
        \State \Comment{Case: \( K = \mathbb{Q}(\sqrt{2})\) }
        \For{\( \ell \in \text{list1} \)}
            \State \( \text{prod} \gets \text{prod} \cdot \ell^{1 + \text{valuation}(4i, \ell)} \)
        \EndFor
    \ElsIf{\( r \) is odd and \( K \subseteq \text{CyclotomicField}(r^2) \)}
        \State \Comment{Case: Odd degree and cyclotomic subfield}
        \For{\( \ell \in \text{list1} \)}
            \State \( \text{prod} \gets \text{prod} \cdot \ell^{1 + \text{valuation}(2ri, \ell)} \)
        \EndFor
    \ElsIf{\( \text{cond} \) is prime and \( K \subseteq \text{CyclotomicField}(\text{cond}) \)}
        \State \Comment{Case: Cyclotomic field with prime conductor}
        \State \( \text{list2} \gets \{ \ell \in \text{list1} \mid \ell \neq 2, \ell \neq \text{cond} \} \)
        \For{\( \ell \in \text{list2} \)}
            \State \( \text{prod} \gets \text{prod} \cdot \ell^{1 + \text{valuation}(i/2, \ell)} \)
        \EndFor
        \If{\( (i \cdot r) \mod (\text{cond} - 1) = 0 \)}
            \State \( \text{prod} \gets \text{prod} \cdot 2^{3 + \text{valuation}(i/2, 2)} \cdot \text{cond}^{1 + \text{valuation}(i/2, \text{cond})} \)
        \Else
            \State \( \text{prod} \gets \text{prod} \cdot 2^{3 + \text{valuation}(i/2, 2)} \)
        \EndIf
    \Else
        \State \Comment{Default case}
        \For{\( \ell \in \text{list1} \)}
            \State \( \text{prod} \gets \text{prod} \cdot \ell^{1 + \text{valuation}(2i, \ell)} \)
        \EndFor
    \EndIf
    \State \Return \( \text{prod} \)
\EndFunction
\end{algorithmic}
\end{algorithm}

\vspace{0.5cm}

\begin{algorithm}[H]
\caption{Compute Zeta Function Value for Number Field with $ Gal(K/\mathbb{Q}) \cong \mathbb{Z}/p\mathbb{Z}$}
\begin{algorithmic}[1]
\State \textbf{Input}: Number field \( K \), conductor \( \text{cond} \), complex number \( s \)
\State \textbf{Output}: Zeta function value
\Function{Zeta}{$K, \text{cond}, s$}
    \State Compute \( n \gets 1 - s \)
    \State Compute degree \( \deg \gets K.\text{degree}() \)
    \State Create Dirichlet group \( G \gets \text{DirichletGroup}(\text{cond}) \)
    \State Initialize \( \text{Lval} \gets 1 \)
    \If{\( \deg = 2 \)}
        \State \Comment{Quadratic field case}
        \State Compute \( \text{Lval} \gets -\text{QuadraticBernoulliNumber}(n, \text{cond})/n \)
    \Else
        \State \Comment{General degree case}
        \State Compute character \( \chi \gets G.0^{(G.\text{order}()/\deg)} \)
        \For{\( i \) in \( [1, \deg-1] \)}
            \State Compute \( \text{Lval} \gets \text{Lval} \cdot (-(\chi^i).\text{bernoulli}(n)/n) \)
        \EndFor
    \EndIf
    \State \Return \( \zeta(s) \cdot \text{Lval} \)
\EndFunction
\end{algorithmic}
\end{algorithm}

\vspace{0.5cm}

\begin{algorithm}[H]
\caption{Compute the size of $K$-groups for Number Field with $ Gal(K/\mathbb{Q}) \cong \mathbb{Z}/p\mathbb{Z}$ }
\begin{algorithmic}[1]
\State \textbf{Input}: Number field \( K \), conductor \( \text{cond} \), integer \( n \)
\State \textbf{Output}: Rounded K-group size
\Function{KSize}{$K, \text{cond}, n$}
    \State Compute degree \( r \gets K.\text{degree}() \)
    \If{\( n \mod 2 = 0 \)}
        \State Compute \( k \gets (n + 2)/4 \)
        \If{\( k \mod 2 = 0 \)}
            \State \Comment{Even \( k \)}
            \State Compute \( SS \gets (1/2^r) \cdot w(K, \text{cond}, 2k) \cdot \text{Zeta}(K, \text{cond}, 1-2k) \)
        \Else
            \State \Comment{Odd \( k \)}
            \State Compute \( SS \gets (-1)^r \cdot w(K, \text{cond}, 2k) \cdot \text{Zeta}(K, \text{cond}, 1-2k) \)
        \EndIf
    \Else
        \State Compute \( k \gets (n + 1)/4 \)
        \If{\( k \mod 2 = 0 \)}
            \State \Comment{Even \( k \)}
            \State Compute \( SS \gets w(K, \text{cond}, 2k) \)
        \Else
            \State \Comment{Odd \( k \)}
            \State Compute \( SS \gets 2^r \cdot w(K, \text{cond}, 2k) \)
        \EndIf
    \EndIf
    \State \Return \( \text{round}(SS) \)
\EndFunction
\end{algorithmic}
\end{algorithm}

\vspace{0.5cm}

\begin{algorithm}[H]
\caption{Compute Scaled Value for $K$-group}
\begin{algorithmic}[1]
\State \textbf{Input}: Number field \( K \), conductor \( \text{cond} \), integer \( n \)
\State \textbf{Output}: Scaled rational value
\Function{OMG}{$K, \text{cond}, n$}
    \State Compute degree \( r \gets \deg(K) \)
    \State Compute \( k \gets \lfloor (n + 2)/4 \rfloor \)
    \If{\( k \) is even}
        \State \Comment{Even case}
        \State Compute \( SS \gets w(K, \text{cond}, 2k) / 2^r \)
    \Else
        \State \Comment{Odd case}
        \State Compute \( SS \gets (-1)^r \cdot w(K, \text{cond}, 2k) \)
    \EndIf
    \State \Return \( \mathbb{Q}(\mathbb{R}(SS)) \)
\EndFunction
\end{algorithmic}
\end{algorithm}

\vspace{0.5cm}

\subsection{\texttt{KSizeppower(K, cond, p, k, n)}}
\begin{algorithm}[H]
\caption{Compute Size of Higher Even $K$-groups}
\begin{algorithmic}[1]
\State \textbf{Input}: Number field \( K \), conductor \( \text{cond} \), prime \( p \), integers \( k \), \( n \)
\State \textbf{Output}: Size of K-group
\Function{KSizeppower}{$K, \text{cond}, p, k, n$}
    \State Compute degree \( r \gets \deg(K) \)
    \State Compute \( a \gets \lfloor (n + 2)/4 \rfloor \)
    \State Create Dirichlet group \( G \gets \text{DirichletGroup}(\text{cond}) \)
    \State Collect \( \text{characters} \gets \{ \chi \in G \mid \chi.\text{multiplicative\_order}() \in \{1, p\} \} \)
    \State Compute subfields \( \text{Subfields} \gets K.\text{subfields}() \)
    \State Initialize empty list \( \text{inter} \)
    \State \Comment{Collect subfields of degree \( p \)}
    \For{each subfield \( S \in \text{Subfields} \)}
        \If{\( \deg(S) = p \)}
            \State Append \( S \) to \( \text{inter} \)
        \EndIf
    \EndFor
    \State Initialize \( \text{prod} \gets 1 \)
    \State \Comment{Compute Zeta value}
    \For{each \( \chi \in \text{characters} \)}
        \State \( \text{prod} \gets \text{prod} \cdot (-\chi.\text{primitive\_character}().\text{bernoulli}(2a)/(2a)) \)
    \EndFor
   
    \State \( \text{prod} \gets \text{prod} \cdot \zeta(1 - 2a)^{|\text{inter}| - 1} \)
    \State Initialize \( ww \gets 1 \)
    \State \Comment{Compute product of $w_j(F)$}
    \For{each \( i \in \text{inter} \)}
        \State \( ww \gets ww \cdot \text{OMG}(i, i.\text{conductor}(), n) \)
    \EndFor
    \State \Comment{Compute denominator}
    \State \( \text{denominator} \gets \text{KZ}(n)^{(p^k - p)/(p - 1)} \)
    \State \Return \( \text{prod} \cdot ww / \text{denominator} \)
\EndFunction
\end{algorithmic}
\end{algorithm}

\newpage 

\section{Size of higher even $K$-groups}\label{AppB}

In this appendix, we give several tables of the values of the $K$-groups. Whenever possible, we will give both the actual order of the $K$-groups and their prime factorization. For orders that are too big, we will only list either their sizes or their prime factorization.

\begin{center}
\begin{longtable}{|l|l|l|}
\caption{$K_{2}(\Op_K)$ for real quadratic field with discriminant $\leq 501$ } \label{tab:K2}\\

\hline \multicolumn{1}{|c|}{$D$} & \multicolumn{1}{c|}{$K$} & \multicolumn{1}{c|}{ $|K_2(\Op_K)|$} \\ \hline 
\endfirsthead
\multicolumn{3}{c}%
{{\bfseries \tablename\ \thetable{} -- continued from previous page}} \\
\hline \multicolumn{1}{|c|}{$D$} & \multicolumn{1}{c|}{$K$} & \multicolumn{1}{c|}{ $|K_2(\Op_K)|$} \\ \hline 
\endhead

\hline \multicolumn{3}{|r|}{{Continued on next page}} \\ \hline
\endfoot

\hline \hline
\endlastfoot
5  & $\mathbb{Q}(\sqrt{ 5 })$ & $ 4 = 2^2 $ \\ \hline
8  & $\mathbb{Q}(\sqrt{ 8 })$ & $ 4 = 2^2 $ \\ \hline
12  & $\mathbb{Q}(\sqrt{ 12 })$ & $ 4 = 2^2 $ \\ \hline
13  & $\mathbb{Q}(\sqrt{ 13 })$ & $ 4 = 2^2 $ \\ \hline
17  & $\mathbb{Q}(\sqrt{ 17 })$ & $ 8 = 2^3 $ \\ \hline
21  & $\mathbb{Q}(\sqrt{ 21 })$ & $ 8 = 2^3 $ \\ \hline
24  & $\mathbb{Q}(\sqrt{ 24 })$ & $ 12 = 2^2 \cdot 3 $ \\ \hline
28  & $\mathbb{Q}(\sqrt{ 28 })$ & $ 16 = 2^4 $ \\
\hline
\end{longtable}
\end{center}

\newpage

\begin{center}
\begin{longtable}{|l|l|l|}
\caption{$K_{6}(\Op_K)$ for real quadratic field with discriminant $\leq 501$ } \label{tab:K6}\\

\hline \multicolumn{1}{|c|}{$D$} & \multicolumn{1}{c|}{$K$} & \multicolumn{1}{c|}{ $|K_6(\Op_K)|$} \\ \hline 
\endfirsthead
\multicolumn{3}{c}%
{{\bfseries \tablename\ \thetable{} -- continued from previous page}} \\
\hline \multicolumn{1}{|c|}{$D$} & \multicolumn{1}{c|}{$K$} & \multicolumn{1}{c|}{ $|K_6(\Op_K)|$} \\ \hline 
\endhead

\hline \multicolumn{3}{|r|}{{Continued on next page}} \\ \hline
\endfoot

\hline \hline
\endlastfoot
5  & $\mathbb{Q}(\sqrt{ 5 })$ & $ 1 = 1 $ \\ \hline
8  & $\mathbb{Q}(\sqrt{ 8 })$ & $ 11 = 11 $ \\ \hline
12  & $\mathbb{Q}(\sqrt{ 12 })$ & $ 23 = 23 $ \\ \hline
13  & $\mathbb{Q}(\sqrt{ 13 })$ & $ 29 = 29 $ \\ \hline
17  & $\mathbb{Q}(\sqrt{ 17 })$ & $ 82 = 2 \cdot 41 $ \\ \hline
21  & $\mathbb{Q}(\sqrt{ 21 })$ & $ 154 = 2 \cdot 7 \cdot 11 $ \\ \hline
24  & $\mathbb{Q}(\sqrt{ 24 })$ & $ 261 = 3^2 \cdot 29 $ \\ \hline
28  & $\mathbb{Q}(\sqrt{ 28 })$ & $ 452 = 2^2 \cdot 113 $ \\
\hline

\end{longtable}
\end{center}

\newpage

\begin{center}
\begin{longtable}{|l|l|l|}
\caption{$K_{10}(\Op_K)$ for real quadratic field with discriminant $\leq 501$ } \label{tab:K10}\\

\hline \multicolumn{1}{|c|}{$D$} & \multicolumn{1}{c|}{$K$} & \multicolumn{1}{c|}{ $|K_{10}(\Op_K)|$} \\ \hline 
\endfirsthead
\multicolumn{3}{c}%
{{\bfseries \tablename\ \thetable{} -- continued from previous page}} \\
\hline \multicolumn{1}{|c|}{$D$} & \multicolumn{1}{c|}{$K$} & \multicolumn{1}{c|}{ $|K_{10}(\Op_K)|$} \\ \hline 
\endhead

\hline \multicolumn{3}{|r|}{{Continued on next page}} \\ \hline
\endfoot

\hline \hline
\endlastfoot
5  & $\mathbb{Q}(\sqrt{ 5 })$ & $ 268 = 2^2 \cdot 67 $ \\ \hline
8  & $\mathbb{Q}(\sqrt{ 8 })$ & $ 1444 = 2^2 \cdot 19^2 $ \\ \hline
12  & $\mathbb{Q}(\sqrt{ 12 })$ & $ 6724 = 2^2 \cdot 41^2 $ \\ \hline
13  & $\mathbb{Q}(\sqrt{ 13 })$ & $ 133852 = 2^2 \cdot 109 \cdot 307 $ \\ \hline
17  & $\mathbb{Q}(\sqrt{ 17 })$ & $ 46328 = 2^3 \cdot 5791 $ \\ \hline
21  & $\mathbb{Q}(\sqrt{ 21 })$ & $ 143768 = 2^3 \cdot 17971 $ \\ \hline
24  & $\mathbb{Q}(\sqrt{ 24 })$ & $ 304332 = 2^2 \cdot 3 \cdot 7 \cdot 3623 $ \\ \hline
28  & $\mathbb{Q}(\sqrt{ 28 })$ & $ 711376 = 2^4 \cdot 173 \cdot 257 $ \\
\hline

\end{longtable}
\end{center}

\newpage

\begin{center}
\begin{longtable}{|l|l|l|}
\caption{$K_{14}(\Op_K)$ for real quadratic field with discriminant $\leq 501$ } \label{tab:K14}\\

\hline \multicolumn{1}{|c|}{$D$} & \multicolumn{1}{c|}{$K$} & \multicolumn{1}{c|}{ $|K_{14}(\Op_K)|$} \\ \hline 
\endfirsthead
\multicolumn{3}{c}%
{{\bfseries \tablename\ \thetable{} -- continued from previous page}} \\
\hline \multicolumn{1}{|c|}{$D$} & \multicolumn{1}{c|}{$K$} & \multicolumn{1}{c|}{ $|K_{14}(\Op_K)|$} \\ \hline 
\endhead

\hline \multicolumn{3}{|r|}{{Continued on next page}} \\ \hline
\endfoot

\hline \hline
\endlastfoot
5  & $\mathbb{Q}(\sqrt{ 5 })$ & $ 361 = 19^2 $ \\ \hline
8  & $\mathbb{Q}(\sqrt{ 8 })$ & $ 24611 = 24611 $ \\ \hline
12  & $\mathbb{Q}(\sqrt{ 12 })$ & $ 257543 = 11 \cdot 13 \cdot 1801 $ \\ \hline
13  & $\mathbb{Q}(\sqrt{ 13 })$ & $ 467669 = 467669 $ \\ \hline
17  & $\mathbb{Q}(\sqrt{ 17 })$ & $ 59901794 = 2 \cdot 19 \cdot 1576363 $ \\ \hline
21  & $\mathbb{Q}(\sqrt{ 21 })$ & $ 17058634 = 2 \cdot 8529317 $ \\ \hline
24  & $\mathbb{Q}(\sqrt{ 24 })$ & $ 46620501 = 3 \cdot 73 \cdot 193 \cdot 1103 $ \\ \hline
28  & $\mathbb{Q}(\sqrt{ 28 })$ & $ 148163732 = 2^2 \cdot 37040933 $ \\
\hline

\end{longtable}
\end{center}

\newpage

\begin{center}
\begin{longtable}{|l|l|l|}
\caption{$K_{18}(\Op_K)$ for real quadratic field with discriminant $\leq 501$ } \label{tab:K18}\\

\hline \multicolumn{1}{|c|}{$D$} & \multicolumn{1}{c|}{$K$} & \multicolumn{1}{c|}{ $|K_{18}(\Op_K)|$} \\ \hline 
\endfirsthead
\multicolumn{3}{c}%
{{\bfseries \tablename\ \thetable{} -- continued from previous page}} \\
\hline \multicolumn{1}{|c|}{$D$} & \multicolumn{1}{c|}{$K$} & \multicolumn{1}{c|}{ $|K_{18}(\Op_K)|$} \\ \hline 
\endhead

\hline \multicolumn{3}{|r|}{{Continued on next page}} \\ \hline
\endfoot

\hline \hline
\endlastfoot

5  & $\mathbb{Q}(\sqrt{ 5 })$ & $ 1651004 = 2^2 \cdot 191 \cdot 2161 $ \\ \hline
8  & $\mathbb{Q}(\sqrt{ 8 })$ & $ 11492164 = 2^2 \cdot 2873041 $ \\ \hline
12  & $\mathbb{Q}(\sqrt{ 12 })$ & $ 270549124 = 2^2 \cdot 83 \cdot 829 \cdot 983 $ \\ \hline
13  & $\mathbb{Q}(\sqrt{ 13 })$ & $ 578188684 = 2^2 \cdot 1289 \cdot 112139 $ \\ \hline
17  & $\mathbb{Q}(\sqrt{ 17 })$ & $ 7408086008 = 2^3 \cdot 23 \cdot 2017 \cdot 19961 $ \\ \hline
21  & $\mathbb{Q}(\sqrt{ 21 })$ & $ 55040678168 = 2^3 \cdot 7 \cdot 982869253 $ \\ \hline
24  & $\mathbb{Q}(\sqrt{ 24 })$ & $ 195898531212 = 2^2 \cdot 3^3 \cdot 23 \cdot 78864143 $ \\ \hline
28  & $\mathbb{Q}(\sqrt{ 28 })$ & $ 847291693456 = 2^4 \cdot 13 \cdot 4073517757 $ \\
\hline

\end{longtable}
\end{center}

\newpage
\begin{landscape}

\begin{center}
\begin{longtable}{|l|l|l|}
\caption{$K_{22}(\Op_K)$ for real quadratic field with discriminant $\leq 501$ } \label{tab:K22}\\

\hline \multicolumn{1}{|c|}{$D$} & \multicolumn{1}{c|}{$K$} & \multicolumn{1}{c|}{ $|K_{22}(\Op_K)|$} \\ \hline 
\endfirsthead
\multicolumn{3}{c}%
{{\bfseries \tablename\ \thetable{} -- continued from previous page}} \\
\hline \multicolumn{1}{|c|}{$D$} & \multicolumn{1}{c|}{$K$} & \multicolumn{1}{c|}{ $|K_{22}(\Op_K)|$} \\ \hline 
\endhead

\hline \multicolumn{3}{|r|}{{Continued on next page}} \\ \hline
\endfoot

\hline \hline
\endlastfoot
5  & $\mathbb{Q}(\sqrt{ 5 })$ & $ 795286411 = 691 \cdot 1150921 $ \\ \hline
8  & $\mathbb{Q}(\sqrt{ 8 })$ & $ 354029435201 = 13^2 \cdot 691 \cdot 3031619 $ \\ \hline
12  & $\mathbb{Q}(\sqrt{ 12 })$ & $ 18752521534133 = 691 \cdot 27138236663 $ \\ \hline
13  & $\mathbb{Q}(\sqrt{ 13 })$ & $ 47067236077919 = 43 \cdot 691 \cdot 1584062063 $ \\ \hline
17  & $\mathbb{Q}(\sqrt{ 17 })$ & $ 1029776943473062 = 2 \cdot 503 \cdot 691 \cdot 23563 \cdot 62869 $ \\ \hline
21  & $\mathbb{Q}(\sqrt{ 21 })$ & $ 11692064980308574 = 2 \cdot 59 \cdot 691 \cdot 8053 \cdot 17806291 $ \\ \hline
24  & $\mathbb{Q}(\sqrt{ 24 })$ & $ 54313104382971231 = 3 \cdot 7 \cdot 349 \cdot 691 \cdot 47869 \cdot 224041 $ \\ \hline
28  & $\mathbb{Q}(\sqrt{ 28 })$ & $ 319737460276467692 = 2^2 \cdot 691 \cdot 1355989 \cdot 85309877 $ \\
\hline

\end{longtable}
\end{center}
\end{landscape}

\newpage\scriptsize	
\begin{landscape}
\begin{center}
\begin{longtable}{|l|l|l|}
\caption{$K_{26}(\Op_K)$ for real quadratic field with discriminant $\leq 501$ } \label{tab:K26}\\

\hline \multicolumn{1}{|c|}{$D$} & \multicolumn{1}{c|}{$K$} & \multicolumn{1}{c|}{ $|K_{26}(\Op_K)|$} \\ \hline 
\endfirsthead
\multicolumn{3}{c}%
{{\bfseries \tablename\ \thetable{} -- continued from previous page}} \\
\hline \multicolumn{1}{|c|}{$D$} & \multicolumn{1}{c|}{$K$} & \multicolumn{1}{c|}{ $|K_{26}(\Op_K)|$} \\ \hline 
\endhead

\hline \multicolumn{3}{|r|}{{Continued on next page}} \\ \hline
\endfoot

\hline \hline
\endlastfoot
5  & $\mathbb{Q}(\sqrt{ 5 })$ & $ 2274367372 =2^2 \cdot 17 \cdot 33446579 $ \\ \hline
8  & $\mathbb{Q}(\sqrt{ 8 })$ & $ 518282899684 = 2^2 \cdot 11 \cdot 151 \cdot 78007661 $ \\ \hline
12  & $\mathbb{Q}(\sqrt{ 12 })$ & $ 61768772694724 = 2^2 \cdot 25321 \cdot 609857161 $ \\ \hline
13  & $\mathbb{Q}(\sqrt{ 13 })$ & $ 181983071666764 = 2^2 \cdot 45495767916691 $ \\ \hline
17  & $\mathbb{Q}(\sqrt{ 17 })$ & $ 6806267386180088 = 2^3 \cdot 547 \cdot 887 \cdot 7411 \cdot 236609 $ \\ \hline
21  & $\mathbb{Q}(\sqrt{ 21 })$ & $ 117965801626742168 = 2^3 \cdot 11 \cdot 31 \cdot 59 \cdot 657911 \cdot 1114019 $ \\ \hline
24  & $\mathbb{Q}(\sqrt{ 24 })$ & $ 715605902169238092 = 2^2 \cdot 3 \cdot 4903 \cdot 6143 \cdot 1979931929 $ \\ \hline
28  & $\mathbb{Q}(\sqrt{ 28 })$ & $ 5733973276382944336 = 2^4 \cdot 23 \cdot 31126933 \cdot 500577719 $ \\
\hline

\end{longtable}
\end{center}
\end{landscape}

\newpage\scriptsize	
\begin{landscape}

\begin{center}
\begin{longtable}{|l|l|l|}
\caption{$K_{30}(\Op_K)$ for real quadratic field with discriminant $\leq 501$ } \label{tab:K30}\\

\hline \multicolumn{1}{|c|}{$D$} & \multicolumn{1}{c|}{$K$} & \multicolumn{1}{c|}{ $|K_{30}(\Op_K)|$} \\ \hline 
\endfirsthead
\multicolumn{3}{c}%
{{\bfseries \tablename\ \thetable{} -- continued from previous page}} \\
\hline \multicolumn{1}{|c|}{$D$} & \multicolumn{1}{c|}{$K$} & \multicolumn{1}{c|}{ $|K_{30}(\Op_K)|$} \\ \hline 
\endhead

\hline \multicolumn{3}{|r|}{{Continued on next page}} \\ \hline
\endfoot

\hline \hline
\endlastfoot
5  & $\mathbb{Q}(\sqrt{ 5 })$ & $ 54701427071177 = 457 \cdot 3617 \cdot 33092833 $ \\ \hline
8  & $\mathbb{Q}(\sqrt{ 8 })$ & $ 159549339299844787 = 79 \cdot 3617 \cdot 558366571709 $ \\ \hline
12  & $\mathbb{Q}(\sqrt{ 12 })$ & $ 42783818174313146311 = 43 \cdot 3617 \cdot 275082254819381 $ \\ \hline
13  & $\mathbb{Q}(\sqrt{ 13 })$ & $ 147939984968716630933 = 47 \cdot 53 \cdot 107 \cdot 3617 \cdot 153454474477 $ \\ \hline
17  & $\mathbb{Q}(\sqrt{ 17 })$ & $ 9460958212613770449794 = 2 \cdot 23 \cdot 109 \cdot 3617 \cdot 19603 \cdot 65599 \cdot 405679 $ \\ \hline
21  & $\mathbb{Q}(\sqrt{ 21 })$ & $ 250243420568581420180298 = 2 \cdot 7 \cdot 37 \cdot 79 \cdot 109 \cdot 379 \cdot 3617 \cdot 40925274407 $ \\ \hline
24  & $\mathbb{Q}(\sqrt{ 24 })$ & $ 1982642779357551544707477 = 3^2 \cdot 3617 \cdot 34157 \cdot 50321 \cdot 87473 \cdot 405089 $ \\ \hline
28  & $\mathbb{Q}(\sqrt{ 28 })$ & $ 21623187890291628882736564 = 2^2 \cdot 3617 \cdot 1494552660374041255373 $ \\
\hline

\end{longtable}
\end{center}
\end{landscape}

\newpage\scriptsize	
\begin{landscape}

\begin{center}
\begin{longtable}{|l|l|l|}
\caption{$K_{34}(\Op_K)$ for real quadratic field with discriminant $\leq 501$ } \label{tab:K34}\\

\hline \multicolumn{1}{|c|}{$D$} & \multicolumn{1}{c|}{$K$} & \multicolumn{1}{c|}{ $|K_{34}(\Op_K)|$} \\ \hline 
\endfirsthead
\multicolumn{3}{c}%
{{\bfseries \tablename\ \thetable{} -- continued from previous page}} \\
\hline \multicolumn{1}{|c|}{$D$} & \multicolumn{1}{c|}{$K$} & \multicolumn{1}{c|}{ $|K_{34}(\Op_K)|$} \\ \hline 
\endhead

\hline \multicolumn{3}{|r|}{{Continued on next page}} \\ \hline
\endfoot

\hline \hline
\endlastfoot

5  & $\mathbb{Q}(\sqrt{ 5 })$ & $ 2285452676758582852 = 2^2 \cdot 41 \cdot 43867 \cdot 317680421579 $ \\ \hline
8  & $\mathbb{Q}(\sqrt{ 8 })$ & $ 3412978697571815263468 = 2^2 \cdot 43867 \cdot 19450718635716001 $ \\ \hline
12  & $\mathbb{Q}(\sqrt{ 12 })$ & $ 2059209895348860985907308 = 2^2 \cdot 11 \cdot 43867 \cdot 1066866320794499171 $ \\ \hline
13  & $\mathbb{Q}(\sqrt{ 13 })$ & $ 108637352378881443421380724 = 2^2 \cdot 1117 \cdot 43867 \cdot 2536679 \cdot 218505594101 $ \\ \hline
17  & $\mathbb{Q}(\sqrt{ 17 })$ & $ 913874565439232931127004456 = 2^3 \cdot 43867 \cdot 2604106063325600483071 $ \\ \hline
21  & $\mathbb{Q}(\sqrt{ 21 })$ & $ 36886276027065642822515422856 = 2^3 \cdot 43867 \cdot 7006659691 \cdot 15001195315081 $ \\ \hline
24  & $\mathbb{Q}(\sqrt{ 24 })$ & $ 381702971297204503144930406724 = 2^2 \cdot 3 \cdot 7 \cdot 1831 \cdot 43867 \cdot 1454939 \cdot 3542771 \cdot 10975697 $ \\ \hline
28  & $\mathbb{Q}(\sqrt{ 28 })$ & $ 5666232102621315072962254685872 = 2^4 \cdot 29 \cdot 21529 \cdot 23801 \cdot 43867 \cdot 543274577837461 $ \\
\hline

\end{longtable}
\end{center}
\end{landscape}

\newpage\scriptsize	
\begin{landscape}

\begin{center}
\begin{longtable}{|l|l|l|}
\caption{$K_{38}(\Op_K)$ for real quadratic field with discriminant $\leq 501$ } \label{tab:K38}\\

\hline \multicolumn{1}{|c|}{$D$} & \multicolumn{1}{c|}{$K$} & \multicolumn{1}{c|}{ $|K_{38}(\Op_K)|$} \\ \hline 
\endfirsthead
\multicolumn{3}{c}%
{{\bfseries \tablename\ \thetable{} -- continued from previous page}} \\
\hline \multicolumn{1}{|c|}{$D$} & \multicolumn{1}{c|}{$K$} & \multicolumn{1}{c|}{ $|K_{38}(\Op_K)|$} \\ \hline 
\endhead

\hline \multicolumn{3}{|r|}{{Continued on next page}} \\ \hline
\endfoot

\hline \hline
\endlastfoot

5  & $\mathbb{Q}(\sqrt{ 5 })$ & $  283 \cdot 617 \cdot 564172514549641 $ \\ \hline
8  & $\mathbb{Q}(\sqrt{ 8 })$ & $  283 \cdot 617 \cdot 211202599 \cdot 51060226589 $ \\ \hline
12  & $\mathbb{Q}(\sqrt{ 12 })$ & $  13 \cdot 47 \cdot 283 \cdot 617 \cdot 23960194641453842173 $ \\ \hline
13  & $\mathbb{Q}(\sqrt{ 13 })$ & $ 31 \cdot 283 \cdot 449 \cdot 617 \cdot 498679 \cdot 10045304252989 $ \\ \hline
17  & $\mathbb{Q}(\sqrt{ 17 })$ & $ 2 \cdot 47 \cdot 283 \cdot 617 \cdot 274787 \cdot 107576263 \cdot 4692574163 $ \\ \hline
21  & $\mathbb{Q}(\sqrt{ 21 })$ & $  2 \cdot 283 \cdot 617 \cdot 32789 \cdot 4191181 \cdot 2921988505217893 $ \\ \hline
24  & $\mathbb{Q}(\sqrt{ 24 })$ & $  3 \cdot 283 \cdot 617 \cdot 3618222022911296408947958807 $ \\ \hline
28  & $\mathbb{Q}(\sqrt{ 28 })$ & $ 2^2 \cdot 283 \cdot 617 \cdot 12391 \cdot 4424992276888190779824023 $ \\
\hline

\end{longtable}
\end{center}
\end{landscape}

\newpage\scriptsize	
\begin{landscape}

\begin{center}
\begin{longtable}{|l|l|}
\caption{$K_{n}(\Op_K)$ for $K=\mathbb{Q}(\sqrt{5})$ } \label{tab:K=Q(sqrt(5))}\\

\hline  \multicolumn{1}{|c|}{$K_{n}(\Op_K)$} & \multicolumn{1}{c|}{ $|K_{n}(\Op_K)|$} \\ \hline 
\endfirsthead
\multicolumn{2}{c}%
{{\bfseries \tablename\ \thetable{} -- continued from previous page}} \\
\hline  \multicolumn{1}{|c|}{$K_{n}(\Op_K)$} & \multicolumn{1}{c|}{ $|K_{n}(\Op_K)|$} \\ \hline 
\endhead

\hline \multicolumn{2}{|r|}{{Continued on next page}} \\ \hline
\endfoot

\hline \hline
\endlastfoot
$ K_{ 2 }(\mathcal{O}_K) $ & $ 4 $ \\ \hline
$ K_{ 6 }(\mathcal{O}_K) $ & $ 1 $ \\ \hline
$ K_{ 10 }(\mathcal{O}_K) $ & $ 268 $ \\ \hline
$ K_{ 14 }(\mathcal{O}_K) $ & $ 361 $ \\ \hline
$ K_{ 18 }(\mathcal{O}_K) $ & $ 1651004 $ \\ \hline
$ K_{ 22 }(\mathcal{O}_K) $ & $ 795286411 $ \\ \hline
$ K_{ 26 }(\mathcal{O}_K) $ & $ 2274367372 $ \\ \hline
$ K_{ 30 }(\mathcal{O}_K) $ & $ 54701427071177 $ \\
\hline

\end{longtable}
\end{center}
\end{landscape}

\newpage
\normalsize

\begin{center} $K_{4k-2}(\Op_E)$ for $E=\Q(\sqrt{2},\sqrt{3},\sqrt{5})$ and $k=1,2,..., 10$ \end{center}  \label{multiquad1}

\begin{tabular}{|c|l|}
    \hline
$|K_{2}(\Op_E)|$ & $2^{11}\cdot 3^2\cdot 7\cdot 17$  \\  \hline
 $|K_{6}(\Op_E)|$ & $2^3\cdot 3^4\cdot 11\cdot 19\cdot 23\cdot 29\cdot 83\cdot 179\cdot 36451$  \\   \hline
 $|K_{10}(\Op_E)|$ & $2^{11}\cdot 3^2\cdot 7\cdot19^2\cdot 41^2\cdot 67\cdot 661\cdot 1303\cdot 1481\cdot 3623\cdot 203999\cdot 1264807$  \\  \hline
$|K_{14}(\Op_E)|$ & $2^3\cdot 3^2\cdot 7\cdot 11\cdot 13\cdot 19^3\cdot 23 \cdot 61\cdot 73\cdot 127\cdot 193\cdot 1103\cdot 1801\cdot 2087 \cdot 24611\cdot 163011019 \cdot$ \\
&$4072124178091$  \\  \hline
$|K_{18}(\Op_E)|$ & $2^{11}\cdot 3^6\cdot 23\cdot 47\cdot 67\cdot 71\cdot 83\cdot 191\cdot 311\cdot 829\cdot 983\cdot 1019\cdot 2161\cdot 2635487\cdot 2873041\cdot  $ \\
& $78864143\cdot 91895159\cdot 106943943050555777$  \\  \hline
$|K_{22}(\Op_E)|$ &
$2^3\cdot 3^2\cdot 7\cdot 13^2\cdot 31\cdot 59\cdot 349\cdot 691\cdot 1201\cdot 47869\cdot 224041\cdot 1150921\cdot 3031619\cdot $  \\
\hline \end{tabular}

\newpage
\normalsize
\begin{center} $K_{4k-2}(\Op_E)$ for $E=\Q(\sqrt{2},\sqrt{3},\sqrt{7})$ and $k=1,2,..., 10$. \end{center}
\[
\begin{tabular}{|c|l|}
  \hline
$|K_{2}(\Op_E)|$ & $2^{13} \cdot 3^4 \cdot 5 $   \\  \hline
 $|K_{6}(\Op_E)|$&  $2^5 \cdot 3^5 \cdot 5^4 \cdot 7^2 \cdot 11^2 \cdot 23 \cdot 29 \cdot 113 \cdot 2503 $\\   \hline
 $|K_{10}(\Op_E)|$& $2^{13} \cdot 3^3 \cdot 5 \cdot 7 \cdot 19^2 \cdot 41^2 \cdot 173 \cdot 257 \cdot 3623 \cdot 17971 \cdot 802729 \cdot 187933043
$ \\  \hline
$|K_{14}(\Op_E)|$ & $2^5 \cdot 3^6 \cdot 5 \cdot 11 \cdot 13 \cdot 23 \cdot 41 \cdot 73 \cdot 193 \cdot 461 \cdot 607 \cdot 1103 \cdot 1801 \cdot $\\&$ 2083 \cdot 24611 \cdot 8529317 \cdot 37040933 \cdot 1019576843
$  \\\hline
$|K_{18}(\Op_E)|$ & $2^{13} \cdot 3^7 \cdot 5 \cdot 7^2 \cdot 13 \cdot 23 \cdot 43 \cdot 83 \cdot 173 \cdot 613 \cdot 727 \cdot 829 \cdot 983 \cdot 1013 \cdot $\\ &$ 1373 \cdot 2873041 \cdot 78864143 \cdot 982869253 \cdot 4073517757 \cdot 15337101636817
$   \\  \hline
$|K_{22}(\Op_E)|$ & $2^5 \cdot 3^3 \cdot 5 \cdot 7 \cdot 13^2 \cdot 59 \cdot 349 \cdot 691 \cdot 7187 \cdot 8053 \cdot 33521 \cdot 47869 \cdot $\\
\hline \end{tabular}
\]

\newpage
\normalsize
\begin{center} $K_{4k-2}(\Op_E)$ for $E=\Q(\sqrt{2},\sqrt{3},\sqrt{11})$ and $k=1,2,..., 10$. \end{center}
\[
\begin{tabular}{|c|l|}
  \hline
$|K_{2}(\Op_E)|$ & $ 2^{13} \cdot 3^2 \cdot 7^2 \cdot 23  $   \\  \hline
$|K_{6}(\Op_E)|$ & $ 2^4 \cdot 3^4 \cdot 11 \cdot 23 \cdot 29^2 \cdot 47 \cdot 2153 \cdot 4967 \cdot 24889  $   \\  \hline
$|K_{10}(\Op_E)|$ & $ 2^{14} \cdot 3^2 \cdot 7 \cdot 11^2 \cdot 19^2 \cdot 41^2 \cdot 277 \cdot 3623 \cdot 31729 \cdot 74231 \cdot 2130727 \cdot 115452991  $   \\  \hline
$|K_{14}(\Op_E)|$ & $ 2^4 \cdot 3^2 \cdot 7^3 \cdot 11 \cdot 13^3 \cdot 23 \cdot 47 \cdot 61 \cdot 73 \cdot 149 \cdot 193 \cdot 211 \cdot 1103 \cdot 1667 \cdot 1801 \cdot 24611 \cdot 30781 \cdot $ \\ & $ 1393361 \cdot 27289513 \cdot 477244787  $   \\  \hline
$|K_{18}(\Op_E)|$ & $ 2^{13} \cdot 3^6 \cdot 23 \cdot 83 \cdot 163 \cdot 829 \cdot 983 \cdot 1151 \cdot 105263 \cdot 568877 \cdot 2873041 \cdot 49103843 \cdot 78864143 \cdot $ \\ & $18701536759 \cdot 487490584799 \cdot 15515203176007  $   \\  \hline
$|K_{22}(\Op_E)|$ & $ 2^4 \cdot 3^2 \cdot 7 \cdot 13^2 \cdot 17 \cdot 73 \cdot 113 \cdot 349 \cdot 691 \cdot 3221 \cdot 8699 \cdot 47869 \cdot 224041 \cdot 1374041 \cdot 3031619 \cdot $ \\
\hline

 \end{tabular}
\]
\newpage
\normalsize
\begin{center} $K_{4k-2}(\Op_E)$ for $E=\Q(\sqrt{2},\sqrt{3},\sqrt{13})$ and $k=1,2,..., 10$. \end{center}
\[
\begin{tabular}{|c|l|}
  \hline

$|K_{2}(\Op_E)|$ & $ 2^{11} \cdot 3^2 \cdot 5^2 \cdot 13 \cdot 23  $   \\  \hline
$|K_{6}(\Op_E)|$ & $ 2^3 \cdot 3^4 \cdot 5^2 \cdot 11 \cdot 23 \cdot 29^2 \cdot 31^3 \cdot 59 \cdot 1409 \cdot 3697  $   \\  \hline
$|K_{10}(\Op_E)|$ & $ 2^{11} \cdot 3^2 \cdot 5 \cdot 7^2 \cdot 13 \cdot 19^2 \cdot 41^3 \cdot 107 \cdot 109 \cdot 307 \cdot 3623 \cdot 34747 \cdot 80382319 \cdot 413993303  $   \\  \hline
$|K_{14}(\Op_E)|$ & $ 2^3 \cdot 3^2 \cdot 5 \cdot 11 \cdot 13 \cdot 73 \cdot 193 \cdot 223 \cdot 1103 \cdot 1481 \cdot 1801 \cdot 24611 \cdot 79367 \cdot 467669 \cdot 35078137 \cdot  $ \\ & $  1967171429 \cdot 7883043022507  $   \\  \hline
$|K_{18}(\Op_E)|$ & $ 2^{11} \cdot 3^6 \cdot 5 \cdot 23 \cdot 83 \cdot 113 \cdot 229 \cdot 233 \cdot 379 \cdot 569 \cdot 829 \cdot 983 \cdot 1129 \cdot 1289 \cdot 2069 \cdot 7573 \cdot 30253 \cdot $ \\ & $  112139 \cdot 975071 \cdot 2873041 \cdot 12322993 \cdot 78864143 \cdot 29455502271467  $   \\  \hline
$|K_{22}(\Op_E)|$ & $ 2^3 \cdot 3^2 \cdot 5 \cdot 7^2 \cdot 13^2 \cdot 37 \cdot 43 \cdot 241 \cdot 349 \cdot 691 \cdot 47869 \cdot 224041 \cdot 3031619 \cdot 5317951 \cdot 1584062063 \cdot $ \\
\hline

 \end{tabular}
\]

\newpage
\normalsize
\begin{center} $K_{4k-2}(\Op_E)$ for $E=\Q(\sqrt{2},\sqrt{3},\sqrt{17})$ and $k=1,2,..., 10$. \end{center}
\[
\begin{tabular}{|c|l|}
  \hline
$|K_{2}(\Op_E)|$ & $ 2^{12} \cdot 3^2 \cdot 13 \cdot 23 \cdot 103  $   \\  \hline
$|K_{6}(\Op_E)|$ & $ 2^4 \cdot 3^4 \cdot 5 \cdot 11 \cdot 23 \cdot 29 \cdot 41 \cdot 409 \cdot 5197 \cdot 6449 \cdot 57241  $   \\  \hline
$|K_{10}(\Op_E)|$ & $ 2^{12} \cdot 3^2 \cdot 7 \cdot 19^2 \cdot 31 \cdot 41^2 \cdot 107 \cdot 167 \cdot 269 \cdot 3623 \cdot 5791 \cdot 159739 \cdot 4956901 \cdot 1640393453  $   \\  \hline
$|K_{14}(\Op_E)|$ & $ 2^4 \cdot 3^2 \cdot 5 \cdot 11 \cdot 13 \cdot 19^2 \cdot 29^2 \cdot 43 \cdot 53 \cdot 73 \cdot 193 \cdot 349 \cdot 1103 \cdot 1801 \cdot 1949 \cdot 24611 \cdot 82279 \cdot  $ \\ & $ 104549 \cdot 1322953 \cdot 1576363 \cdot 2387843 \cdot 8861147  $   \\  \hline
$|K_{18}(\Op_E)|$ & $ 2^{12} \cdot 3^6 \cdot 23^2 \cdot 37 \cdot 61 \cdot 67 \cdot 83 \cdot 167 \cdot 563 \cdot 829 \cdot 983 \cdot 2017 \cdot 19961 \cdot 220903 \cdot 252383 \cdot 2873041 \cdot  $ \\ & $ 78864143 \cdot 1685953303 \cdot 116420774341 \cdot 16552191488784001  $   \\  \hline
$|K_{22}(\Op_E)|$ & $ 2^4 \cdot 3^2 \cdot 5^2 \cdot 7 \cdot 13^2 \cdot 179 \cdot 349 \cdot 503 \cdot 523 \cdot 691 \cdot 3697 \cdot 23563 \cdot 33629 \cdot 47869 \cdot 62869 \cdot 224041 \cdot  $ \\
\hline

 \end{tabular}
\]
\newpage
\normalsize
\begin{center} $K_{4k-2}(\Op_E)$ for $E=\Q(\sqrt{2},\sqrt{3},\sqrt{19})$ and $k=1,2,..., 10$. \end{center}
\[
\begin{tabular}{|c|l|}
  \hline
$|K_{2}(\Op_E)|$ & $ 2^{12} \cdot 3^2 \cdot 7 \cdot 11 \cdot 19 \cdot 41  $   \\  \hline
$|K_{6}(\Op_E)|$ & $ 2^5 \cdot 3^4 \cdot 5^2 \cdot 11 \cdot 19 \cdot 23^2 \cdot 29 \cdot 47 \cdot 61 \cdot 109 \cdot 137 \cdot 571 \cdot 1429  $   \\  \hline
$|K_{10}(\Op_E)|$ & $ 2^{12} \cdot 3^2 \cdot 7 \cdot 19^2 \cdot 29 \cdot 41^2 \cdot 373 \cdot 1327 \cdot 3391 \cdot 3623 \cdot 43067 \cdot 794483 \cdot 1488671 \cdot 5222837  $   \\  \hline
$|K_{14}(\Op_E)|$ & $ 2^7 \cdot 3^2 \cdot 5^3 \cdot 7 \cdot 11^2 \cdot 13 \cdot 73 \cdot 79^2 \cdot 193 \cdot 257 \cdot 1103 \cdot 1801 \cdot 3061 \cdot 24611 \cdot 21865969 \cdot  $ \\ & $ 2195956421 \cdot 3353773067 \cdot 24276385561  $   \\  \hline
$|K_{18}(\Op_E)|$ & $ 2^{12} \cdot 3^6 \cdot 13 \cdot 19^2 \cdot 23^2 \cdot 41 \cdot 59 \cdot 61^2 \cdot 79 \cdot 83 \cdot 223 \cdot 677 \cdot 829 \cdot 983 \cdot 6823 \cdot 44647 \cdot 46279 \cdot $ \\ & $  2873041 \cdot 7718891 \cdot 12186413 \cdot 78864143 \cdot 15176499719 \cdot 700276785671  $   \\  \hline
$|K_{22}(\Op_E)|$ & $ 2^5 \cdot 3^2 \cdot 5^2 \cdot 7 \cdot 11^2 \cdot 13^2 \cdot 17 \cdot 151 \cdot 229 \cdot 271 \cdot 349 \cdot 691 \cdot 787 \cdot 6367 \cdot 7559 \cdot 47869 \cdot 224041 \cdot $ \\
\hline

 \end{tabular}
\]

\newpage

\bt Let $K$ be a cyclic cubic field with conductor $f$, where $f$ is congruent to 1 modulo 3 . Then there exist a unique integer $a \equiv 2(\bmod 3)$ for each $f$ such that $K=\mathbb{Q}(\theta)$ where $\theta$ is a root of
$$
P(X)=X^3-3 f X-f a
$$

More generally, Hasse indicated that one can write $f=\left(a^2+3 b^2\right) / 4$, where $a$ and $b$ satisfy the conditions
$$
\begin{gathered}
a \equiv 2(\bmod 3), b \equiv 0(\bmod 3), b>0 \text { if } 3 \nmid f \\
a \equiv 6(\bmod 9), b \equiv 3 \operatorname{or} 6(\bmod 9), b>0 \text { if } 3 \mid f .
\end{gathered}
$$

\et

\begin{center}
\begin{longtable}{|l|l|l|l|}
\caption{$K_{2}(\Op_K)$ for cyclic cubic field with conductor $\leq 500$ } \label{tab:cubicK2}\\

\hline \multicolumn{1}{|c|}{$f$} & \multicolumn{1}{c|}{$a$} & \multicolumn{1}{c|}{$b$} & \multicolumn{1}{c|}{ $|K_2(\Op_K)|$} \\ \hline 
\endfirsthead
\multicolumn{4}{c}%
{{\bfseries \tablename\ \thetable{} -- continued from previous page}} \\
\hline \multicolumn{1}{|c|}{$f$} & \multicolumn{1}{c|}{$a$} & \multicolumn{1}{c|}{$b$} & \multicolumn{1}{c|}{ $|K_2(\Op_K)|$} \\ \hline 
\endhead

\hline \multicolumn{4}{|r|}{{Continued on next page}} \\ \hline
\endfoot

\hline \hline
\endlastfoot
7 & -1 & 3 & $ 2^3 $ \\ \hline
9 & -3 & 3 & $ 2^3 $ \\ \hline
13 & 5 & 3 & $ 2^3 $ \\ \hline
19 & -7 & 3 & $ 2^3 \cdot 3 $ \\ \hline
31 & -4 & 6 & $ 2^5 \cdot 7 $ \\ \hline
37 & 11 & 3 & $ 2^3 \cdot 3 \cdot 7 $ \\ \hline
43 & 8 & 6 & $ 2^5 \cdot 19 $ \\ \hline
61 & -1 & 9 & $ 2^3 \cdot 7 \cdot 19 $ \\
\hline

\end{longtable}
\end{center}

\newpage
\begin{center}
\begin{longtable}{|l|l|l|l|}
\caption{$K_{6}(\Op_K)$ for cyclic cubic field with conductor $\leq 500$ } \label{tab:cubicK6}\\

\hline \multicolumn{1}{|c|}{$f$} & \multicolumn{1}{c|}{$a$} & \multicolumn{1}{c|}{$b$} & \multicolumn{1}{c|}{ $|K_6(\Op_K)|$} \\ \hline 
\endfirsthead
\multicolumn{4}{c}%
{{\bfseries \tablename\ \thetable{} -- continued from previous page}} \\
\hline \multicolumn{1}{|c|}{$f$} & \multicolumn{1}{c|}{$a$} & \multicolumn{1}{c|}{$b$} & \multicolumn{1}{c|}{ $|K_6(\Op_K)|$} \\ \hline 
\endhead

\hline \multicolumn{4}{|r|}{{Continued on next page}} \\ \hline
\endfoot

\hline \hline
\endlastfoot
7 & -1 & 3 & $ 79 $ \\ \hline
9 & -3 & 3 & $ 199 $ \\ \hline
13 & 5 & 3 & $ 103 \cdot 109 $ \\ \hline
19 & -7 & 3 & $ 3 \cdot 61 \cdot 67 $ \\ \hline
31 & -4 & 6 & $ 2^2 \cdot 151 \cdot 757 $ \\ \hline
37 & 11 & 3 & $ 3 \cdot 433513 $ \\ \hline
43 & 8 & 6 & $ 2^2 \cdot 421 \cdot 2683 $ \\ \hline
61 & -1 & 9 & $ 7 \cdot 13 \cdot 397 \cdot 1237 $ \\
\hline

\end{longtable}
\end{center}

\newpage
\begin{center}
\begin{longtable}{|l|l|l|l|}
\caption{$K_{10}(\Op_K)$ for cyclic cubic field with conductor $\leq 500$ } \label{tab:cubicK10}\\

\hline \multicolumn{1}{|c|}{$f$} & \multicolumn{1}{c|}{$a$} & \multicolumn{1}{c|}{$b$} & \multicolumn{1}{c|}{ $|K_{10}(\Op_K)|$} \\ \hline 
\endfirsthead
\multicolumn{4}{c}%
{{\bfseries \tablename\ \thetable{} -- continued from previous page}} \\
\hline \multicolumn{1}{|c|}{$f$} & \multicolumn{1}{c|}{$a$} & \multicolumn{1}{c|}{$b$} & \multicolumn{1}{c|}{ $|K_{10}(\Op_K)|$} \\ \hline 
\endhead

\hline \multicolumn{4}{|r|}{{Continued on next page}} \\ \hline
\endfoot

\hline \hline
\endlastfoot

7 & -1 & 3 & $ 2^3 \cdot 7393 $ \\ \hline
9 & -3 & 3 & $ 2^3 \cdot 7 \cdot 43 \cdot 1171 $ \\ \hline
13 & 5 & 3 & $ 2^3 \cdot 73 \cdot 91807 $ \\ \hline
19 & -7 & 3 & $ 2^3 \cdot 3 \cdot 127 \cdot 21720427 $ \\ \hline
31 & -4 & 6 & $ 2^5 \cdot 13 \cdot 1915270939 $ \\ \hline
37 & 11 & 3 & $ 2^3 \cdot 3 \cdot 109 \cdot 2034277813 $ \\ \hline
43 & 8 & 6 & $ 2^5 \cdot 910644093019 $ \\ \hline
61 & -1 & 9 & $ 2^3 \cdot 163 \cdot 443347 \cdot 2260603 $ \\
\hline

\end{longtable}
\end{center}

\newpage
\begin{center}
\begin{longtable}{|l|l|l|l|}
\caption{$K_{14}(\Op_K)$ for cyclic cubic field with conductor $\leq 500$ } \label{tab:cubicK14}\\

\hline \multicolumn{1}{|c|}{$f$} & \multicolumn{1}{c|}{$a$} & \multicolumn{1}{c|}{$b$} & \multicolumn{1}{c|}{ $|K_{14}(\Op_K)|$} \\ \hline 
\endfirsthead
\multicolumn{4}{c}%
{{\bfseries \tablename\ \thetable{} -- continued from previous page}} \\
\hline \multicolumn{1}{|c|}{$f$} & \multicolumn{1}{c|}{$a$} & \multicolumn{1}{c|}{$b$} & \multicolumn{1}{c|}{ $|K_{14}(\Op_K)|$} \\ \hline 
\endhead

\hline \multicolumn{4}{|r|}{{Continued on next page}} \\ \hline
\endfoot

\hline \hline
\endlastfoot

7 & -1 & 3 & $ 142490119 $ \\ \hline
9 & -3 & 3 & $ 2648750959 $ \\ \hline
13 & 5 & 3 & $ 19 \cdot 109 \cdot 1377630997 $ \\ \hline
19 & -7 & 3 & $ 3 \cdot 13 \cdot 37 \cdot 571 \cdot 78998317 $ \\ \hline
31 & -4 & 6 & $ 2^2 \cdot 7 \cdot 19 \cdot 109 \cdot 3109 \cdot 564653359 $ \\ \hline
37 & 11 & 3 & $ 3 \cdot 7 \cdot 37 \cdot 2251 \cdot 374977 \cdot 2179921 $ \\ \hline
43 & 8 & 6 & $ 2^2 \cdot 3445737657466136863 $ \\ \hline
61 & -1 & 9 & $ 7 \cdot 369202062414708268357 $ \\
\hline

\end{longtable}
\end{center}

\newpage
\begin{center}
\begin{longtable}{|l|l|l|l|}
\caption{$K_{18}(\Op_K)$ for cyclic cubic field with conductor $\leq 500$ } \label{tab:cubicK18}\\

\hline \multicolumn{1}{|c|}{$f$} & \multicolumn{1}{c|}{$a$} & \multicolumn{1}{c|}{$b$} & \multicolumn{1}{c|}{ $|K_{18}(\Op_K)|$} \\ \hline 
\endfirsthead
\multicolumn{4}{c}%
{{\bfseries \tablename\ \thetable{} -- continued from previous page}} \\
\hline \multicolumn{1}{|c|}{$f$} & \multicolumn{1}{c|}{$a$} & \multicolumn{1}{c|}{$b$} & \multicolumn{1}{c|}{ $|K_{18}(\Op_K)|$} \\ \hline 
\endhead

\hline \multicolumn{4}{|r|}{{Continued on next page}} \\ \hline
\endfoot

\hline \hline
\endlastfoot

7 & -1 & 3 & $ 2^3 \cdot 1141452324871 $ \\ \hline
9 & -3 & 3 & $ 2^3 \cdot 619 \cdot 93657511669 $ \\ \hline
13 & 5 & 3 & $ 2^3 \cdot 79 \cdot 1051 \cdot 251878440139 $ \\ \hline
19 & -7 & 3 & $ 2^3 \cdot 3 \cdot 79 \cdot 164186413 \cdot 727327213 $ \\ \hline
31 & -4 & 6 & $ 2^5 \cdot 61 \cdot 541 \cdot 73015779967328830657 $ \\ \hline
37 & 11 & 3 & $ 2^3 \cdot 3 \cdot 67 \cdot 3187 \cdot 13955699568289894873 $ \\ \hline
43 & 8 & 6 & $ 2^5 \cdot 14142451 \cdot 443576761 \cdot 6210365209 $ \\ \hline
61 & -1 & 9 & $ 2^3 \cdot 7^2 \cdot 43 \cdot 21475957 \cdot 502049533 \cdot 5253179929 $ \\
\hline

\end{longtable}
\end{center}

\newpage
\begin{center}
\begin{longtable}{|l|l|l|l|}
\caption{$K_{22}(\Op_K)$ for cyclic cubic field with conductor $\leq 500$ } \label{tab:cubicK22}\\

\hline \multicolumn{1}{|c|}{$f$} & \multicolumn{1}{c|}{$a$} & \multicolumn{1}{c|}{$b$} & \multicolumn{1}{c|}{ $|K_{22}(\Op_K)|$} \\ \hline 
\endfirsthead
\multicolumn{4}{c}%
{{\bfseries \tablename\ \thetable{} -- continued from previous page}} \\
\hline \multicolumn{1}{|c|}{$f$} & \multicolumn{1}{c|}{$a$} & \multicolumn{1}{c|}{$b$} & \multicolumn{1}{c|}{ $|K_{22}(\Op_K)|$} \\ \hline 
\endhead

\hline \multicolumn{4}{|r|}{{Continued on next page}} \\ \hline
\endfoot

\hline \hline
\endlastfoot

7 & -1 & 3 & $ 691 \cdot 10903 \cdot 278995143079 $ \\ \hline
9 & -3 & 3 & $ 7 \cdot 199 \cdot 691 \cdot 727 \cdot 139801 \cdot 20873329 $ \\ \hline
13 & 5 & 3 & $ 691 \cdot 4640714412917661501919 $ \\ \hline
19 & -7 & 3 & $ 3 \cdot 127 \cdot 691 \cdot 1051 \cdot 1579 \cdot 89983 \cdot 9569822168437 $ \\ \hline
31 & -4 & 6 & $ 2^2 \cdot 691 \cdot 556499562270997220800171338307 $ \\ \hline
37 & 11 & 3 & $ 3 \cdot 19 \cdot 691 \cdot 22937767 \cdot 3683860021705494153552817 $ \\ \hline
43 & 8 & 6 & $ 2^2 \cdot 31 \cdot 691 \cdot 322259933988079 \cdot 103359961664870167 $ \\ \hline
61 & -1 & 9 & $ 13 \cdot 109 \cdot 691 \cdot 6427 \cdot 1409607628208664285027824452561 $ \\
\hline

\end{longtable}
\end{center}

\newpage
\begin{landscape}
\begin{center}
\begin{longtable}{|l|l|l|l|}
\caption{$K_{26}(\Op_K)$ for cyclic cubic field with conductor $\leq 500$ } \label{tab:cubicK26}\\

\hline \multicolumn{1}{|c|}{$f$} & \multicolumn{1}{c|}{$a$} & \multicolumn{1}{c|}{$b$} & \multicolumn{1}{c|}{ $|K_{26}(\Op_K)|$} \\ \hline 
\endfirsthead
\multicolumn{4}{c}%
{{\bfseries \tablename\ \thetable{} -- continued from previous page}} \\
\hline \multicolumn{1}{|c|}{$f$} & \multicolumn{1}{c|}{$a$} & \multicolumn{1}{c|}{$b$} & \multicolumn{1}{c|}{ $|K_{26}(\Op_K)|$} \\ \hline 
\endhead

\hline \multicolumn{4}{|c|}{{Continued on next page}} \\ \hline
\endfoot

\hline \hline
\endlastfoot
7 & -1 & 3 & $ 2^3 \cdot 1033 \cdot 5410539334962035689 $ \\ \hline
9 & -3 & 3 & $ 2^3 \cdot 6229 \cdot 6703 \cdot 7252358598014773 $ \\ \hline
13 & 5 & 3 & $ 2^3 \cdot 433 \cdot 4780577201281149751995847 $ \\ \hline
19 & -7 & 3 & $ 2^3 \cdot 3 \cdot 9631 \cdot 241921 \cdot 886987 \cdot 9407223976677403 $ \\ \hline
31 & -4 & 6 & $ 2^5 \cdot 7 \cdot 35221988143 \cdot 2662419387637 \cdot 12220658621851 $ \\ \hline
37 & 11 & 3 & $ 2^3 \cdot 3 \cdot 7 \cdot 73 \cdot 17633023 \cdot 1014439300939 \cdot 138947032097099251 $ \\ \hline
43 & 8 & 6 & $ 2^5 \cdot 151 \cdot 811 \cdot 19348005730989956984586749257774654477 $ \\ \hline
61 & -1 & 9 & $ 2^3 \cdot 7 \cdot 73 \cdot 691 \cdot 7861484998413330473769625373284270201903 $ \\
\hline

\end{longtable}
\end{center}
\end{landscape}

\newpage
\begin{landscape}
\begin{center}
\begin{longtable}{|l|l|l|l|}
\caption{$K_{30}(\Op_K)$ for cyclic cubic field with conductor $\leq 500$ } \label{tab:cubicK30}\\

\hline \multicolumn{1}{|c|}{$f$} & \multicolumn{1}{c|}{$a$} & \multicolumn{1}{c|}{$b$} & \multicolumn{1}{c|}{ $|K_{30}(\Op_K)|$} \\ \hline 
\endfirsthead
\multicolumn{4}{c}%
{{\bfseries \tablename\ \thetable{} -- continued from previous page}} \\
\hline \multicolumn{1}{|c|}{$f$} & \multicolumn{1}{c|}{$a$} & \multicolumn{1}{c|}{$b$} & \multicolumn{1}{c|}{ $|K_{30}(\Op_K)|$} \\ \hline 
\endhead

\hline \multicolumn{4}{|r|}{{Continued on next page}} \\ \hline
\endfoot

\hline \hline
\endlastfoot
7 & -1 & 3 & $ 3617 \cdot 19387 \cdot 6997171 \cdot 399890401961287 $ \\ \hline
9 & -3 & 3 & $ 307 \cdot 3617 \cdot 60133 \cdot 3045261617686098608209 $ \\ \hline
13 & 5 & 3 & $ 79 \cdot 727 \cdot 3617 \cdot 6379 \cdot 20023 \cdot 2964700454864089513687 $ \\ \hline
19 & -7 & 3 & $ 3 \cdot 2749 \cdot 3617 \cdot 216133 \cdot 120665009424320663984357967631 $ \\ \hline
31 & -4 & 6 & $ 2^2 \cdot 61 \cdot 757 \cdot 3617 \cdot 11779 \cdot 385367476974136869940539620801617129 $ \\ \hline
37 & 11 & 3 & $ 3 \cdot 3617 \cdot 1331645369667008529811 \cdot 50580184115960583692794603 $ \\ \hline
43 & 8 & 6 & $ 2^2 \cdot 1021 \cdot 3617 \cdot 3548773 \cdot 21467863 \cdot 68521369216236798313780551013657 $ \\ \hline
61 & -1 & 9 & $ 7 \cdot 13 \cdot 3617 \cdot 7879 \cdot 4427083 \cdot 93484663 \cdot 3665145242262922690883456381675179 $ \\
\hline

\end{longtable}
\end{center}
\end{landscape}

\newpage
\begin{landscape}
\begin{center}
\begin{longtable}{|l|l|l|l|}
\caption{$K_{34}(\Op_K)$ for cyclic cubic field with conductor $\leq 500$ } \label{tab:cubicK34}\\

\hline \multicolumn{1}{|c|}{$f$} & \multicolumn{1}{c|}{$a$} & \multicolumn{1}{c|}{$b$} & \multicolumn{1}{c|}{ $|K_{34}(\Op_K)|$} \\ \hline 
\endfirsthead
\multicolumn{4}{c}%
{{\bfseries \tablename\ \thetable{} -- continued from previous page}} \\
\hline \multicolumn{1}{|c|}{$f$} & \multicolumn{1}{c|}{$a$} & \multicolumn{1}{c|}{$b$} & \multicolumn{1}{c|}{ $|K_{34}(\Op_K)|$} \\ \hline 
\endhead

\hline \multicolumn{4}{|r|}{{Continued on next page}} \\ \hline
\endfoot

\hline \hline
\endlastfoot

7 & -1 & 3 & $ 2^3 \cdot 97 \cdot 43867 \cdot 9105835027474306843301627809 $ \\ \hline
9 & -3 & 3 & $ 2^3 \cdot 7 \cdot 14173 \cdot 17377 \cdot 43867 \cdot 10156277756426731543364053 $ \\ \hline
13 & 5 & 3 & $ 2^3 \cdot 31 \cdot 457 \cdot 43867 \cdot 160103028995840019789566403113309713 $ \\ \hline
19 & -7 & 3 & $ 2^3 \cdot 3 \cdot 127 \cdot 43867 \cdot 3492296987054368270811659441517033240521273 $ \\ \hline
31 & -4 & 6 & $ 2^5 \cdot 13 \cdot 211 \cdot 43867 \cdot 31314686098652351761 \cdot 106977022094660181288313338049 $ \\ \hline
37 & 11 & 3 & $ 2^3 \cdot 3 \cdot 1471 \cdot 43867 \cdot 158796450281985779346247 \cdot 25653447464008652181704392441 $ \\ \hline
43 & 8 & 6 & $ 2^5 \cdot 12613 \cdot 43867 \cdot 2586433380116077 \cdot 1852336684128189769 \cdot 14314013275848210811 $ \\ \hline
61 & -1 & 9 & $ 2^3 \cdot 13873 \cdot 43867 \cdot 122149 \cdot 366703 \cdot 138009121 \cdot 8335310557652067473727098925357618225913 $ \\
\hline
\end{longtable}
\end{center}
\end{landscape}

\newpage
\begin{landscape}
\begin{center}
\begin{longtable}{|l|l|l|l|}
\caption{$K_{38}(\Op_K)$ for cyclic cubic field with conductor $\leq 500$ } \label{tab:cubicK38}\\

\hline \multicolumn{1}{|c|}{$f$} & \multicolumn{1}{c|}{$a$} & \multicolumn{1}{c|}{$b$} & \multicolumn{1}{c|}{ $|K_{38}(\Op_K)|$} \\ \hline 
\endfirsthead
\multicolumn{4}{c}%
{{\bfseries \tablename\ \thetable{} -- continued from previous page}} \\
\hline \multicolumn{1}{|c|}{$f$} & \multicolumn{1}{c|}{$a$} & \multicolumn{1}{c|}{$b$} & \multicolumn{1}{c|}{ $|K_{38}(\Op_K)|$} \\ \hline 
\endhead

\hline \multicolumn{4}{|r|}{{Continued on next page}} \\ \hline
\endfoot

\hline \hline
\endlastfoot

7 & -1 & 3 & $ 283 \cdot 617 \cdot 21766351 \cdot 51183510123014870096951001289 $ \\ \hline
9 & -3 & 3 & $ 283 \cdot 617 \cdot 8260033 \cdot 1043737853428019551309070510936383 $ \\ \hline
13 & 5 & 3 & $ 283 \cdot 617 \cdot 63201616464333885249960435172680450449283126067 $ \\ \hline
19 & -7 & 3 & $ 3 \cdot 13 \cdot 283 \cdot 617 \cdot 42449294191 \cdot 44815116170694531517 \cdot 175398115851660773737 $ \\ \hline
31 & -4 & 6 & $ 2^2 \cdot 7^2 \cdot 283 \cdot 617 \cdot 10986244111 \cdot 239780029764485047 \cdot 152948338622224269898755373835749 $ \\ \hline
37 & 11 & 3 & $ 3 \cdot 7 \cdot 37 \cdot 283 \cdot 617 \cdot 1201 \cdot 13009 \cdot 3624641239 \cdot 14287726884889 \cdot 4021844039046472557853263229573 $ \\ \hline
43 & 8 & 6 & $ 2^2 \cdot 19 \cdot 283 \cdot 617 \cdot 11680233144769270254641089874436388213929946771087386287928112717 $ \\ \hline
61 & -1 & 9 & $ 7 \cdot 19 \cdot 43 \cdot 283 \cdot 617 \cdot 6571 \cdot 27799 \cdot 6054661 \cdot 13895513551$ \\
\hline

\end{longtable}
\end{center}
\end{landscape}

\newpage

\begin{center} $K_{4k-2}(\Op_K)$ for $K=\Q(x^9 - 15x^7 - 4x^6 + 54x^5 + 12x^4 - 38x^3 - 9x^2 + 6x + 1)$. \end{center}
\[
\begin{tabular}{|c|l|}
  \hline
$|K_{2}(\Op_K)|$ & $2^{11} \cdot 7 \cdot 19 \cdot 67$  \\  \hline
 $|K_{6}(\Op_K)|$ & $2^2 \cdot 7 \cdot 61 \cdot 79 \cdot 199 \cdot 38821 \cdot 54922771$  \\   \hline
 $|K_{10}(\Op_K)|$ & $2^{11} \cdot 7 \cdot 37 \cdot 43 \cdot 79 \cdot 1171 \cdot 7393 \cdot 85093 \cdot 715561519 \cdot 79494942841
$  \\  \hline
$|K_{14}(\Op_K)|$ & $2^2 \cdot 7 \cdot 61 \cdot 379 \cdot 523 \cdot 7369 \cdot 7933 \cdot 142490119 \cdot 2648750959 \cdot 11054056927  \cdot 81261959882835397$  \\  \hline 
$|K_{18}(\Op_K)|$ & $2^{11} \cdot 7 \cdot 43 \cdot 397 \cdot 619 \cdot 853  \cdot 93657511669 \cdot 1141452324871 \cdot 541853290098721179727 $ \\& $ \cdot 220279917248393405756637993283
$ \\   \hline
$|K_{22}(\Op_K)|$ &  $2^2 \cdot 7 \cdot 43 \cdot 199 \cdot 337 \cdot 691 \cdot 727 \cdot 1693 \cdot 10903 \cdot 139801 \cdot 20873329 \cdot 176035957 \cdot 278995143079 $ \\& $ \cdot 22631898140730217577524567 \cdot 1860470714192517073619292462626641
$ \\
\hline
\end{tabular}
\]

\newpage
\begin{center} $K_{4k-2}(\Op_K)$ for $K=\Q(x^9 - 20 x^7 - 8 x^6 + 103 x^5 + 46 x^4 - 157 x^3 - 32 x^2 + 60 x + 8)$. \end{center}
\[
\begin{tabular}{|c|l|}
  \hline
$|K_{2}(\Op_K)|$ & $2^{11} \cdot 3^2 \cdot 61 \cdot 151
$  \\  \hline
 $|K_{6}(\Op_K)|$ & $2^2 \cdot 3^2 \cdot 19 \cdot 37 \cdot 79 \cdot 103 \cdot 109 \cdot 348583 \cdot 71601037$  \\   \hline
 $|K_{10}(\Op_K)|$ & $2^{11} \cdot 3^2 \cdot 31 \cdot 73 \cdot 7393 \cdot 8887 \cdot 91807 \cdot 1683637 \cdot 687518473 \cdot 16097131033
$  \\  \hline
$|K_{14}(\Op_K)|$ & $2^2 \cdot 3^2 \cdot 19 \cdot 109 \cdot 16595923 \cdot 23441179 \cdot 95692321 \cdot 142490119 \cdot 225826141 \cdot 1377630997 $\\
& $\cdot 3631091984690689
$  \\  \hline 
$|K_{18}(\Op_K)|$ & $2^{11} \cdot 3^2 \cdot 79 \cdot 109 \cdot 1051 \cdot 5827 \cdot 92041 \cdot 108059683 \cdot 251878440139 \cdot 1141452324871$ \\ & $\cdot 18378228102982243 \cdot 13637640022514423708067259723
$ \\   \hline
$|K_{22}(\Op_K)|$ &  $2^2 \cdot 3^2 \cdot 19 \cdot 571 \cdot 691 \cdot 10903 \cdot 1992259 \cdot 20344123 \cdot 278995143079 \cdot 3984523386229999$\\
\hline
\end{tabular}
\]

\newpage

\begin{center} $K_{4k-2}(\Op_K)$ for $K=\Q(x^9 - 3 x^8 - 18 x^7 + 38 x^6 + 93 x^5 - 147 x^4 - 161 x^3 + 201 x^2 + 57 x - 53)$. \end{center}
\[
\begin{tabular}{|c|l|}
  \hline
$|K_{2}(\Op_K)|$ & $2^{11} \cdot 5^2 \cdot 31 \cdot 433
$  \\  \hline
 $|K_{6}(\Op_K)|$ & $2^2 \cdot 7 \cdot 37 \cdot 103 \cdot 109 \cdot 139 \cdot 199 \cdot 810319 \cdot 180400267
$  \\   \hline
 $|K_{10}(\Op_K)|$ & $2^{11} \cdot 5^2 \cdot 7^2 \cdot 19^2 \cdot 43 \cdot 73 \cdot 1171 \cdot 39703 \cdot 91807 \cdot 29048281 \cdot 289954369 \cdot 550057153
$  \\  \hline
$|K_{14}(\Op_K)|$ & $2^2 \cdot 19 \cdot 37 \cdot 109 \cdot 409 \cdot 4153 \cdot 114021877 \cdot 1377630997 \cdot 2648750959 \cdot 245130870023113$ \\ & $\cdot 294107925568998608917
$  \\  \hline
$|K_{18}(\Op_K)|$ & $2^{11} \cdot 5^2 \cdot 7 \cdot 19 \cdot 61 \cdot 79 \cdot 241 \cdot 619 \cdot 1051 \cdot 2617 \cdot 3069607 \cdot 4598359 \cdot 316641931 \cdot 950769277$ \\ & $\cdot 93657511669 \cdot 251878440139 \cdot 368139927619417921605333439$ \\   \hline
$|K_{22}(\Op_K)|$ & $1604335832984073148836447161751607044808354582986951286912673046681783649472918493 $
\\
\hline
\end{tabular}
\]

\newpage
\begin{center} $K_{4k-2}(\Op_K)$ for $K=\Q(x^9 - 26x^7 - 14x^6 + 181x^5 + 154x^4 - 295x^3 - 224x^2 + 84x + 56)$. \end{center}
\[
\begin{tabular}{|c|l|}
  \hline
$|K_{2}(\Op_K)|$ & $2^{11} \cdot 3^3 \cdot 211 \cdot 463
$  \\  \hline
 $|K_{6}(\Op_K)|$ & $2^2 \cdot 3^3 \cdot 37^2 \cdot 61 \cdot 67 \cdot 79 \cdot 2083 \cdot 13297 \cdot 94328149
$  \\   \hline
 $|K_{10}(\Op_K)|$ & $2^{11} \cdot 3^3 \cdot 61 \cdot 127 \cdot 2593 \cdot 6673 \cdot 7393 \cdot 21720427 \cdot 184874857207 \cdot 111166719688441
 $  \\  \hline
$|K_{14}(\Op_K)|$ & $2^2 \cdot 3^3 \cdot 13 \cdot 37 \cdot 571 \cdot 78998317 \cdot 79761169 \cdot 142490119 \cdot 1292041763952027217 $\\ & $\cdot 26056263648538002464811367
$  \\  \hline
$|K_{18}(\Op_K)|$ & $2^{11} \cdot 3^3 \cdot 13 \cdot 31 \cdot 79 \cdot 619 \cdot 481171 \cdot 1891987 \cdot 8563771 \cdot 164186413 \cdot 727327213  \cdot 2216441803$ \\& $\cdot 93886491697 \cdot 132251313361 \cdot 1141452324871 \cdot 54176532538747 
$ \\   \hline
$|K_{22}(\Op_K)|$ & $2^2 \cdot 3^3 \cdot 127 \cdot 499 \cdot 691 \cdot 1051 \cdot 1579 \cdot 10903 \cdot 89983 \cdot 445447 \cdot 10187959 \cdot 16459123 \cdot 24409453$ \\
\hline
\end{tabular}
\]

\newpage
\begin{center} $K_{4k-2}(\Op_K)$ for $K=\Q(x^9 - 3x^8 - 24x^7 + 74x^6 + 117x^5 - 429x^4 + 119x^3 + 243x^2 - 69x - 37)$. \end{center}
\[
\begin{tabular}{|c|l|}
  \hline
$|K_{2}(\Op_K)|$ & $2^{17} \cdot 3^3 \cdot 7 \cdot 853
$  \\  \hline
 $|K_{6}(\Op_K)|$ & $2^6 \cdot 3^3 \cdot 37 \cdot 61 \cdot 67 \cdot 163 \cdot 199 \cdot 10127737 \cdot 121832929
$  \\   \hline
 $|K_{10}(\Op_K)|$ & $2^{17} \cdot 3^3 \cdot 7 \cdot 43 \cdot 109 \cdot 127 \cdot 199 \cdot 1171 \cdot 1609 \cdot 21720427 \cdot 11611579651003 \cdot 210365363298073
 $  \\  \hline
$|K_{14}(\Op_K)|$ & $2^6 \cdot 3^3 \cdot 7 \cdot 13 \cdot 37 \cdot 571 \cdot 2389 \cdot 16843 \cdot 53097529 \cdot 78998317 \cdot 2648750959 \cdot 2091202347589$ \\& $\cdot 10090761648726379294303087 $  \\  \hline
$|K_{18}(\Op_K)|$ & $2^{17} \cdot 3^3 \cdot 37 \cdot 79 \cdot 619 \cdot 2689 \cdot 2853079 \cdot 164186413 \cdot 727327213 \cdot 93657511669 \cdot 8667944615377 $\\& $\cdot 2018888113095739 \cdot 128090458406494487710527386353591
$  \\   \hline
$|K_{22}(\Op_K)|$ & $2^6 \cdot 3^3 \cdot 7 \cdot 127 \cdot 199 \cdot 691 \cdot 727 \cdot 1051 \cdot 1579 \cdot 89983 \cdot 139801 \cdot 20873329 \cdot 61814887 \cdot 9569822168437 
$\\
\hline
\end{tabular}
\]

\newpage
\begin{center} $K_{4k-2}(\Op_K)$ for $K=\Q(x^9 - 38x^7 - 13x^6 + 409x^5 + 257x^4 - 1333x^3 - 832x^2 + 1143x + 433)$. \end{center}
\[
\begin{tabular}{|c|l|}
  \hline
$|K_{2}(\Op_K)|$ & $2^{11} \cdot 3^3 \cdot 7 \cdot 13 \cdot 37^2 \cdot 61
$  \\  \hline
 $|K_{6}(\Op_K)|$ & $2^2 \cdot 3^4 \cdot 7 \cdot 73^2 \cdot 79 \cdot 151 \cdot 757 \cdot 8923 \cdot 165709 \cdot 22507651
$  \\   \hline
$|K_{10}(\Op_K)|$ & $2^{11} \cdot 3^3 \cdot 7^2 \cdot 13 \cdot 619 \cdot 811 \cdot 7393 \cdot 80713 \cdot 497839 \cdot 1915270939 \cdot 3178658347 \cdot 417756469213$  \\  \hline
$|K_{14}(\Op_K)|$ & $2^2 \cdot 3^3 \cdot 7 \cdot 19 \cdot 37 \cdot 109^2 \cdot 1297 \cdot 3109 \cdot 165811 \cdot 2873701 \cdot 44755957 \cdot 142490119 \cdot 564653359$ \\& $\cdot 196925567948917 \cdot 384845428942442407$ \\   \hline
$|K_{18}(\Op_K)|$ & $2^{11} \cdot 3^4 \cdot 7 \cdot 61 \cdot 241 \cdot 541 \cdot 1747 \cdot 922303 \cdot 1141452324871 \cdot 1521171751711 $ \\& $ \cdot 73015779967328830657 \cdot 443263229448661087922611 \cdot 84116162599285849483190455717
$ \\   \hline
$|K_{22}(\Op_K)|$ & $2^2 \cdot 3^3 \cdot 7 \cdot 691 \cdot 10903 \cdot 680803 \cdot 278995143079 \cdot 556499562270997220800171338307 $ \\
\hline

 \end{tabular}
\]

\newpage
\begin{center} $K_{4k-2}(\Op_K)$ for $K=\Q(x^9 - 32x^7 - 22x^6 + 259x^5 + 140x^4 - 787x^3 - 88x^2 + 804x - 248)$. \end{center}
\[
\begin{tabular}{|c|l|}
  \hline
$|K_{2}(\Op_K)|$ & $2^{11} \cdot 3^3 \cdot 19 \cdot 41^2 \cdot 109
$  \\  \hline
 $|K_{6}(\Op_K)|$ & $2^2 \cdot 3^3 \cdot 13 \cdot 43 \cdot 61 \cdot 67 \cdot 103 \cdot 109 \cdot 733 \cdot 1801 \cdot 254383 \cdot 109999789
$  \\   \hline
 $|K_{10}(\Op_K)|$ & $2^{11} \cdot 3^3 \cdot 73 \cdot 127 \cdot 283 \cdot 19423 \cdot 91807 \cdot 430819 \cdot 21720427 \cdot 8183094133 \cdot 919384681715200627 $  \\  \hline
 $|K_{14}(\Op_K)|$ & $2^2 \cdot 3^3 \cdot 13 \cdot 19 \cdot 31 \cdot 37 \cdot 43 \cdot 109 \cdot 307 \cdot 499 \cdot 547 \cdot 571 \cdot 28591 \cdot 132637 \cdot 12197323 \cdot 78998317 $ \\& $ \cdot 87963529 \cdot 1377630997 \cdot 2854627417 \cdot 240585279948271
$ \\   \hline
 $|K_{18}(\Op_K)|$ & $2^{11} \cdot 3^3 \cdot 61 \cdot 79^2 \cdot 1051 \cdot 8179 \cdot 164186413 \cdot 727327213 \cdot 251878440139 $ \\& $ \cdot 27643166473351105893646084798110199 \cdot 3458222037941198950681890148399088850481
$  \\   \hline
$|K_{22}(\Op_K)|$ &  $2^2 \cdot 3^3 \cdot 13 \cdot 127 \cdot 691 \cdot 811 \cdot 1051 \cdot 1579 \cdot 14071 \cdot 89983 \cdot 19366423 \cdot 69368983 \cdot 2256178447 $ \\
\hline
\end{tabular}\]

\newpage
\begin{center} $K_{4k-2}(\Op_K)$ for $K=\Q(x^9 - 44x^7 - 46x^6 + 559x^5 + 1004x^4 - 1711x^3 - 4600x^2 - 2172x + 232)$. \end{center}
\[
\begin{tabular}{|c|l|}
  \hline
$|K_{2}(\Op_K)|$ & $2^{11} \cdot 3^3 \cdot 7 \cdot 13 \cdot 151 \cdot 2731
$  \\  \hline
 $|K_{6}(\Op_K)|$ & $2^2 \cdot 3^3 \cdot 7 \cdot 79 \cdot 97 \cdot 193 \cdot 433513 \cdot 857929 \cdot 358827911989
$  \\   \hline
 $|K_{10}(\Op_K)|$ & $2^{11} \cdot 3^3 \cdot 7 \cdot 109 \cdot 241 \cdot 7393 \cdot 23539 \cdot 51511 \cdot 158527 \cdot 325921 \cdot 516405157 \cdot 927161959 \cdot 2034277813$  \\  \hline
 $|K_{14}(\Op_K)|$ & $2^2 \cdot 3^3 \cdot 7 \cdot 37 \cdot 1663 \cdot 2251 \cdot 374977 \cdot 2179921 \cdot 7711441 \cdot 10671433 \cdot 18639421 \cdot 27586543  $ \\& $\cdot 142490119 \cdot 144398941 \cdot 127491038349082784491$ \\   \hline
 $|K_{18}(\Op_K)|$ & $2^{11} \cdot 3^3 \cdot 7^3 \cdot 13 \cdot 67 \cdot 79 \cdot 151 \cdot 229 \cdot 313 \cdot 3187 \cdot 7993 \cdot 34849 \cdot 1141452324871 \cdot 4250969512807  $ \\& $ \cdot 15487348934881 \cdot 13955699568289894873 \cdot 4138006003106445147574373965452871$ \\   \hline
$|K_{22}(\Op_K)|$ & $2^2 \cdot 3^3 \cdot 7 \cdot 19 \cdot 61 \cdot 691 \cdot 1447 \cdot 2341 \cdot 10903 \cdot 83857 \cdot 8755189 \cdot 22937767 \cdot 226709731  $ \\
\hline
\end{tabular}
\]

\newpage
\begin{center} $K_{4k-2}(\Op_K)$ for $K=\Q(x^9 - 3x^8 - 36x^7 + 101x^6 + 327x^5 - 921x^4 - 575x^3 + 1632x^2 + 480x - 449)$. \end{center}
\[
\begin{tabular}{|c|l|}
  \hline
$|K_{2}(\Op_K)|$ & $2^2 \cdot 7 \cdot 9973 \cdot 11941
$  \\  \hline
 $|K_{6}(\Op_K)|$ & $2^2 \cdot 19^2 \cdot 151 \cdot 199 \cdot 757 \cdot 37441 \cdot 135043 \cdot 1836171281947
$  \\   \hline
 $|K_{10}(\Op_K)|$ & $2^{11} \cdot 7 \cdot 13 \cdot 37 \cdot 43^2 \cdot 1171 \cdot 3637 \cdot 6763657 \cdot 1915270939 \cdot 441462524953003 \cdot 515902866728209 $  \\  \hline
 $|K_{14}(\Op_K)|$ &  $2^2 \cdot 7 \cdot 19^2 \cdot 109 \cdot 3109 \cdot 364627 \cdot 939931 \cdot 5917897 \cdot 459269077 \cdot 564653359 \cdot 2648750959  $ \\& $\cdot 10214959351483 \cdot 2374239246468122056381$\\   \hline
 $|K_{18}(\Op_K)|$ & $2^{11} \cdot 61 \cdot 139 \cdot 541 \cdot 619 \cdot 24109 \cdot 26959 \cdot 2073199 \cdot 93657511669 \cdot 1188768945769  $ \\& $\cdot 6302924564642056459 \cdot 73015779967328830657 \cdot 124965518382890301922663887688295293$ \\   \hline
$|K_{22}(\Op_K)|$ & $
17666920266552871504913887731620487255344911639572864038875898364322885665747059051 $ \\
\hline\end{tabular}
\]

\newpage
\begin{center} $K_{4k-2}(\Op_K)$ for $K=\Q(x^9 - 50x^7 - 39x^6 + 733x^5 + 999x^4 - 3097x^3 - 5616x^2 + 279x + 2619)$. \end{center}
\[
\begin{tabular}{|c|l|}
  \hline
$|K_{2}(\Op_K)|$ & $2^2 \cdot 3^2 \cdot 19 \cdot 67 \cdot 79 \cdot 4297
$  \\  \hline
 $|K_{6}(\Op_K)|$ & $2^2 \cdot 3^2 \cdot 79 \cdot 421 \cdot 853 \cdot 2683 \cdot 154339 \cdot 6876217 \cdot 1204501843
$  \\   \hline
 $|K_{10}(\Op_K)|$ & $2^{11} \cdot 3^2 \cdot 11^2 \cdot 13^2 \cdot 73 \cdot 769 \cdot 7393 \cdot 90787 \cdot 124567 \cdot 16481791 \cdot 25465987 \cdot 966637069  $ \\& $\cdot 910644093019$  \\  \hline
 $|K_{14}(\Op_K)|$ & $2^2 \cdot 3^2 \cdot 97 \cdot 19482739 \cdot 26954449 \cdot 142490119 \cdot 10119573943 \cdot 79031993305297  $ \\& $\cdot 3445737657466136863 \cdot 11412193713048369162721
$ \\   \hline
 $|K_{18}(\Op_K)|$ & $2^{11} \cdot 3^2 \cdot 31 \cdot 619 \cdot 1597 \cdot 14142451 \cdot 443576761 \cdot 6210365209 \cdot 1141452324871 \cdot 6532269157058137  $ \\& $\cdot 146010655759125116338789 \cdot 11925994175930471072587615508928217129$  \\
\hline\end{tabular}
\]

\newpage

\begin{center} $K_{4k-2}(\Op_K)$ for $K=\Q(x^{27} - 9x^{26} - 45x^{25} + 558x^{24} + 564x^{23} - 14325x^{22} + 3278x^{21} + 200187x^{20} - 152250x^{19} - 1689790x^{18} + 1662003x^{17} + 9058716x^{16} - 9244050x^{15} - 31587348x^{14} + 28970679x^{13} + 71864499x^{12} - 50585529x^{11} - 103930188x^{10} + 44676536x^9 + 88142967x^8 - 15697422x^7 - 37059379x^6 + 2104236x^5 + 7157724x^4 - 134797x^3 - 528300x^2 + 16716x + 7064)$. \end{center}
\[
\begin{tabular}{|c|l|}
  \hline
$|K_{2}(\Op_K)|$ & $2^{35} \cdot 3^{10} \cdot 5^2 \cdot 7^3 \cdot 19 \cdot 31^2 \cdot 61 \cdot 67^2 \cdot 139 \cdot 151 \cdot 211 \cdot 433 \cdot 1753 \cdot 96697$  \\  \hline
 $|K_{6}(\Op_K)|$ & $2^8 \cdot 3^9 \cdot 7^3 \cdot 19^2 \cdot 31 \cdot 37^2 \cdot 43 \cdot 61 \cdot 79 \cdot 103 \cdot 109 \cdot 139 \cdot 181 \cdot 199^2 \cdot 283$ \\& $\cdot 27367 \cdot 38821 \cdot 348583 \cdot 810319 \cdot 23326789 \cdot 54922771 \cdot 71601037$ \\ & $\cdot 180400267 \cdot 719882071 \cdot 142149038377 \cdot 164087854132681$  \\   \hline
 $|K_{10}(\Op_K)|$ & $2^{35} \cdot 3^8 \cdot 5^2 \cdot 7^3 \cdot 19^4 \cdot 31 \cdot 37 \cdot 43 \cdot 73 \cdot 79 \cdot 97 \cdot 163 \cdot 1171 \cdot 7393 \cdot 8887$ \\& $\cdot 10099 \cdot 39703 \cdot 85093 \cdot 91807 \cdot 303421 \cdot 1683637 \cdot 29048281 \cdot 160969201$ \\ & $\cdot 289954369 \cdot 550057153 \cdot 687518473 \cdot 715561519 \cdot 16097131033 \cdot 79494942841$ \\ & $\cdot 518975951743 \cdot 546870623491 \cdot 467812489885669 \cdot 6471445414111759$ \\
\hline
 \end{tabular}%
\]

\newpage

\section{Examples}\label{AppC}

This appendix introduces three methods for computing the size of higher even $K$-groups, each tailored to specific types of number fields. The major distinction among these methods lies in how the special values of the Dedekind zeta function are computed. While software like SageMath and Magma provide built-in functions for numerical approximations of these values, we aim to compute them exactly as fractions. The first method is designed for fields whose degrees are $p$-powers, leveraging their specific arithmetic structures. The second method applies to fields of higher degree or those whose degrees are not $p$-powers, where the associated Dirichlet character must be manually input from the LMFDB to ensure precise fractional computation. Finally, for the quadratic case, we provide a simplified approach using Siegel's formula and also offer a faster method to calculate $w_j(F)$, enhancing efficiency while maintaining exactness. Together, these methods form a comprehensive framework for rigorously determining $K$-group sizes across various field types.

The first method is the most general and theoretically works for every totally real abelian number field. The primary difference from the first method lies in the computation of the Dedekind zeta function, where the associated Dirichlet character must be manually input from the LMFDB to obtain the special values. While this generality ensures broad applicability, the computation of $w_j(F)$ can become significantly slower as the degree of the field increases, due to the amplified complexity and computational effort required. To illustrate this method, we include an example of computing the Dedekind zeta value for a degree 21 field, highlighting both its flexibility and the challenges posed by high-degree fields.
\begin{lstlisting}[language=Python,  captionpos=t, caption= Computing special value of Dedekind Zeta function for $\mathbb{Q}(x^{21} - 7x^{20} - 70x^{19} + 462x^{18} + 2135x^{17} - 12411x^{16} - 36610x^{15} + 175044x^{14} + 373940x^{13} - 1403661x^{12} - 2218069x^{11} + 6566007x^{10} + 6982234x^9 - 17907827x^8 - 9448729x^7 + 26548844x^6 + 686581x^5 - 16732429x^4 + 5281647x^3 + 1717044x^2 - 573440x - 71167)$]
from sage.modular.dirichlet import DirichletCharacter

H = DirichletGroup(637, base_ring=CyclotomicField(2))

M = H._module

chi1 = DirichletCharacter(H, M([0,0]))

from sage.modular.dirichlet import DirichletCharacter

H = DirichletGroup(637, base_ring=CyclotomicField(42))

M = H._module

chi2 = DirichletCharacter(H, M([6,28]))

from sage.modular.dirichlet import DirichletCharacter

H = DirichletGroup(637, base_ring=CyclotomicField(14))

M = H._module

chi3 = DirichletCharacter(H, M([10,0]))

from sage.modular.dirichlet import DirichletCharacter

H = DirichletGroup(637, base_ring=CyclotomicField(6))

M = H._module

chi4 = DirichletCharacter(H, M([0,2]))

from sage.modular.dirichlet import DirichletCharacter

H = DirichletGroup(637, base_ring=CyclotomicField(42))

M = H._module

chi5 = DirichletCharacter(H, M([36,28]))

from sage.modular.dirichlet import DirichletCharacter

H = DirichletGroup(637, base_ring=CyclotomicField(42))

M = H._module

chi6 = DirichletCharacter(H, M([24,14]))

from sage.modular.dirichlet import DirichletCharacter

H = DirichletGroup(637, base_ring=CyclotomicField(14))

M = H._module

chi7 = DirichletCharacter(H, M([6,0]))

from sage.modular.dirichlet import DirichletCharacter

H = DirichletGroup(637, base_ring=CyclotomicField(42))

M = H._module

chi8 = DirichletCharacter(H, M([30,14]))

from sage.modular.dirichlet import DirichletCharacter

H = DirichletGroup(637, base_ring=CyclotomicField(42))

M = H._module

chi9 = DirichletCharacter(H, M([24,28]))

from sage.modular.dirichlet import DirichletCharacter

H = DirichletGroup(637, base_ring=CyclotomicField(14))

M = H._module

chi10 = DirichletCharacter(H, M([2,0]))

from sage.modular.dirichlet import DirichletCharacter

H = DirichletGroup(637, base_ring=CyclotomicField(42))

M = H._module

chi11 = DirichletCharacter(H, M([12,28]))

from sage.modular.dirichlet import DirichletCharacter

H = DirichletGroup(637, base_ring=CyclotomicField(14))

M = H._module

chi12 = DirichletCharacter(H, M([12,0]))

from sage.modular.dirichlet import DirichletCharacter

H = DirichletGroup(637, base_ring=CyclotomicField(42))

M = H._module

chi13 = DirichletCharacter(H, M([6,14]))

from sage.modular.dirichlet import DirichletCharacter

H = DirichletGroup(637, base_ring=CyclotomicField(6))

M = H._module

chi14 = DirichletCharacter(H, M([0,4]))

from sage.modular.dirichlet import DirichletCharacter

H = DirichletGroup(637, base_ring=CyclotomicField(14))

M = H._module

chi15 = DirichletCharacter(H, M([8,0]))

from sage.modular.dirichlet import DirichletCharacter

H = DirichletGroup(637, base_ring=CyclotomicField(42))

M = H._module

chi16 = DirichletCharacter(H, M([36,14]))

from sage.modular.dirichlet import DirichletCharacter

H = DirichletGroup(637, base_ring=CyclotomicField(42))

M = H._module

chi17 = DirichletCharacter(H, M([18,14]))

from sage.modular.dirichlet import DirichletCharacter

H = DirichletGroup(637, base_ring=CyclotomicField(42))

M = H._module

chi18 = DirichletCharacter(H, M([30,28]))

from sage.modular.dirichlet import DirichletCharacter

H = DirichletGroup(637, base_ring=CyclotomicField(14))

M = H._module

chi19 = DirichletCharacter(H, M([4,0]))

from sage.modular.dirichlet import DirichletCharacter

H = DirichletGroup(637, base_ring=CyclotomicField(42))

M = H._module

chi20 = DirichletCharacter(H, M([18,28]))

from sage.modular.dirichlet import DirichletCharacter

H = DirichletGroup(637, base_ring=CyclotomicField(42))

M = H._module

chi21 = DirichletCharacter(H, M([12,14]))




CHI = [chi1,chi2,chi3, chi4, chi5, chi6, chi7,chi8,chi9, chi10, chi11, chi12, chi13, chi14, chi15, chi16, chi17,chi18,chi19, chi20, chi21]

CChi=[]
for i in CHI:
  CChi.append(i.primitive_character())

prod =1

for i in CChi:
  prod*=-i.bernoulli(2)/2
\end{lstlisting}\label{deg21}

For quadratic fields, we employ Siegel’s formula to compute the special values of the Dedekind zeta function. Furthermore, we introduce a faster method to compute $w_j(F)$, optimizing the process while maintaining the accuracy required to determine the sizes of the $K$ groups. This specialized approach significantly improves computational efficiency in the quadratic case.
\begin{lstlisting}[language=Python,  captionpos=t, caption= Computing special value of Dedekind Zeta function using Siegel's formula]
# Define a power series ring over rational numbers
R.<q> = QQ[['q']]

# Compute Eisenstein series q-expansions with given normalization and precision
F4 = eisenstein_series_qexp(4, 4000, normalization='constant')  # Eisenstein series of weight 4
F6 = eisenstein_series_qexp(6, 4000, normalization='constant')  # Eisenstein series of weight 6

# Define the discriminant modular form Delta as a combination of F4 and F6
Delta = (F4^3 - F6^2) / 1728

def T(h):
    """
    Computes a modular form T(h) based on Delta and Eisenstein series of various weights.
    :param h: Integer, typically the weight of the modular form.
    :return: A q-expansion of the modular form T(h).
    """
    if Mod(h, 12) == 2:
        r = floor(h / 12)
    else:
        r = floor(h / 12) + 1

    k = 12 * r - h + 2  # Compute the required weight for the Eisenstein series
    if k == 0:
        T = Delta^(-r)  # Case where k = 0, use only powers of Delta
    else:
        G = eisenstein_series_qexp(k, r + 1, normalization='constant')  # Eisenstein series of weight k
        T = G * Delta^(-r)  # Combine with Delta

    return T + O(q)  # Return the truncated q-expansion

def b(l, h):
    """
    Compute the normalized coefficient b_j(h).
    :param l: Index of the coefficient to retrieve.
    :param h: Weight parameter for T(h).
    :return: thew value of b_j(h).
    """
    C = T(h).coefficients()  # Get coefficients of T(h)
    L = len(C)  # Total number of coefficients
    return -C[L - l - 1] / C[L - 1]  # Return normalized coefficient

def eSiegel(D, j):
    """
    Computes the Siegel-type divisor sum for a given discriminant D and exponent j.
    :param D: Discriminant.
    :param j: Exponent for the divisor function.
    :return: Siegel sum value.
    """
    haha = 0
    bla = 0

    # Compute first term where D/4 is an integer
    if D / 4 in ZZ:
        for i in divisors(round(D / 4)):
            haha += i^j

    # Compute second term for values b where (D - b^2)/4 is an integer
    for b in range(1, floor(sqrt(D)) + 1):
        if (D - b^2) / 4 in ZZ:
            for i in divisors(round((D - b^2) / 4)):
                bla += i^j

    return haha + 2 * bla

def S(F, l, m):
    """
    Computes a sum involving Kronecker symbols and Siegel sums for a field F.
    :param F: A quadratic field.
    :param l: Integer parameter.
    :param m: Integer parameter.
    :return: The computed sum.
    """
    sum = 0
    D = F.discriminant()  # Discriminant of the quadratic field

    # Iterate over divisors of l
    for j in divisors(l):
        sum += kronecker_symbol(D, j) * j^(2 * m - 1) * eSiegel((l / j)^2 * D, 2 * m - 1)

    return sum

def Zeta(D, s):
    """
    Computes a zeta-like function for a totally real quadratic field with discriminant D.
    :param D: Discriminant of the quadratic field.
    :param s: An integer parameter.
    :return: Zeta value.
    """
    F = QuadraticField(D)  # Create the quadratic field
    r = 2  # Degree of the field (always 2 for quadratic fields)
    k = (1 - s) / 2  # Parameter derived from s

    # Determine the range for summation based on k*r modulo 6
    if Mod(k * r, 6) == 1:
        c = floor(k * r / 6)
    else:
        c = floor(k * r / 6) + 1

    sum = 0

    # Compute Dedekind zeta function as in Theorem 4.1
    for j in range(1, c + 1):
        sum += b(j, 2 * k * r) * S(F, j, k)

    return 2^r * sum  # Multiply by 2^r and return the result




def w(K, i):
    """
    Computes w_j(F) for real quadratic fields.
    :param K: A number field, typically quadratic.
    :param i: An integer parameter used to determine specific properties of the primes.
    :return: w_j(F).
    """
    list1 = []

    # Collect primes satisfying (i / (ell - 1)) in Z
    for ell in list(primes(2, i + 2)):
        if (i / (ell - 1)) in ZZ:
            list1.append(ell)

    # Remove 2 from list1 to handle it separately
    list1 = [ell for ell in list1 if ell != 2]

    # Determine additional primes based on k = 2 * i and specific divisibility criteria
    k = 2 * i
    list2 = []
    for ell in list(primes(2, 2 * i + 2)):
        if (k / (ell - 1)) in ZZ and (i / (ell - 1)) not in ZZ:
            list2.append(ell)

    # Compute the product based on the properties of list1 and list2
    prod = 1
    for ell in list1:
        if (i / ell) not in ZZ:
            prod *= ell  # Multiply by prime
        else:
            prod *= ell**(i.valuation(ell) + 1)  # Adjust power based on valuation

    # Handle the contribution from 2 separately
    prod *= 2**(i.valuation(2) + 2)

    # Additional adjustment for primes in list2
    for ell in list2:
        if K == QuadraticField(ell):  # Check if K matches a specific quadratic field
            prod *= ell**(i.valuation(ell) + 1)

    # Adjust for special cases where K is QuadraticField(2) or QuadraticField(8)
    if K == QuadraticField(2) or K == QuadraticField(8):
        prod = 2 * prod

    return prod


def KSize(K, n):
    """
    Computes the size of higher even K-groups for the field K and parameter n.
    :param K: A number field.
    :param n: An integer parameter related to the K-group.
    :return: The computed size as an integer.
    """
    r = K.degree()  # Degree of the field
    D = K.discriminant()  # Discriminant of the field

    # Case when n is even
    if n % 2 == 0:
        k = (n + 2) // 4
        if k % 2 == 0:  # k is even
            SS = (1 / (2^r)) * w(K, 2 * k) * Zeta(D, 1 - 2 * k)
        else:  # k is odd
            SS = (-1)^r * w(K, 2 * k) * Zeta(D, 1 - 2 * k)

    # Case when n is odd
    else:
        k = (n + 1) // 4
        if k % 2 == 0:  # k is even
            SS = w(K, 2 * k)
        else:  # k is odd
            SS = 2^r * w(K, 2 * k)

    return round(SS)  # Return the rounded result
\end{lstlisting}\label{siegel}

\bigskip
For $p$-elementary abelian totally real fields, our code leverages the specific arithmetic properties of these fields to compute the size of higher even $K$-groups efficiently, directly utilizing the structured relationship between the field and its Dedekind zeta function.
\begin{lstlisting}[language=Python, captionpos=t, caption= Code for computing size of higher even $K$-groups]

def is_subfield(K, L):
    """
    Checks if K is a subfield of L by verifying the existence of embeddings.
    :param K: A number field or similar algebraic structure.
    :param L: A number field into which embeddings are checked.
    :return: True if K embeds into L, otherwise False.
    """
    return len(K.embeddings(L)) > 0


#Computes w_i(K) based on Propsition 4.2
def w(K, cond, i):
    """
    Computes a weighted product of prime powers for a given number field and conditions.
    :param K: A number field.
    :param cond: An integer condition, e.g., conductor of a field.
    :param i: An integer parameter used in the calculations.
    :return: A product based on certain arithmetic properties of K and i.
    """
    r = K.degree()  # Degree of the number field K.
    list1 = []

    # Collecting primes that satisfy specific congruence properties.
    for ell in prime_range(2, i + 2):
        if i % (ell - 1) == 0:
            list1.append(ell)

    prod = 1

    # Case 1: K is a subfield of QuadraticField(2) and QuadraticField(8)
    if is_subfield(QuadraticField(2), K) and is_subfield(K, QuadraticField(8)):
        for ell in list1:
            prod *= ell ** (1 + (4 * i).valuation(ell))

    # Case 2: K is related to a CyclotomicField of degree r^2
    elif r % 2 == 1 and is_subfield(K, CyclotomicField(r**2)):
        for ell in list1:
            prod *= ell ** (1 + (2 * r * i).valuation(ell))

    # Case 3: K is related to a CyclotomicField of a prime conductor
    elif is_subfield(K, CyclotomicField(cond)) and is_prime(cond):
        list2 = [ell for ell in list1 if ell != 2 and ell != cond]
        for ell in list2:
            prod *= ell ** (1 + (i / 2).valuation(ell))
        if (i * r) % (cond - 1) == 0:
            prod = 2 ** (3 + (i / 2).valuation(2)) * cond ** (1 + (i / 2).valuation(cond)) * prod
        else:
            prod = 2 ** (3 + (i / 2).valuation(2)) * prod

    # Default case
    else:
        for ell in list1:
            prod *= ell ** (1 + (2 * i).valuation(ell))

    return prod


#The following function computes the special value of the Dedekind Zeta function for number field K.
#But in this case, it only deals with the fields with p degrees. 
#We also have a general version which deals with any fields. (see appendix xxx)
def Zeta(K, cond, s):
    """
    Computes a modified zeta function for the number field K.
    :param K: A number field.
    :param cond: An integer condition, typically related to the field's conductor.
    :param s: A negative integer.
    :return: The special value of Dedekind Zeta fucntion of K.
    """
    n = 1 - s
    deg = K.degree()  # Degree of the field.
    G = DirichletGroup(cond)
    Lval = 1

    if deg == 2:
        Lval = -QuadraticBernoulliNumber(n, cond) / n
    else:
        chi = G.0^(ZZ(G.order() / deg))
        for i in (1..deg - 1):
            Lval *= (-(chi^i).bernoulli(n) / n)

    return zeta(s) * Lval



#Computes the size of even K group of ring of integers of the field K using Theorem 2.3 and Corollary 2.4.
def KSize(K, cond, n):
    """
    Computes the size of even K group of ring of integers of the field K.
    :param K: A number field K.
    :param cond: An integer condition, typically a conductor.
    :param n: An integer parameter.
    :return: The computed size.
    """
    r = K.degree()

    if n % 2 == 0:
        k = (n + 2) // 4
        if k % 2 == 0:
            SS = 1 / (2**r) * w(K, cond, 2 * k) * Zeta(K, cond, 1 - 2 * k)
        else:
            SS = (-1)**r * w(K, cond, 2 * k) * Zeta(K, cond, 1 - 2 * k)
    else:
        k = (n + 1) // 4
        if k % 2 == 0:
            SS = w(K, cond, 2 * k)
        else:
            SS = 2**r * w(K, cond, 2 * k)

    return round(SS)


#This function computes size of even K groups of Z, which appeared in Theorem 4.1. 
def KZ(n):
    """
    Computes size of even K groups of Z, which appeared in Theorem 4.1. 
    :param n: An integer.
    :return: A scaling factor based on Bernoulli numbers and prime valuations.
    """
    ZZ = IntegerRing()

    if n % 8 == 2:
        k = ZZ((n - 2) / 8)
        s = 2 * k + 1
        c = bernoulli(2 * s) * w(QQ, 1, 2 * s) / (4 * s)
        return abs(2 * c)
    elif n % 8 == 6:
        k = ZZ((n - 6) / 8)
        s = 2 * k + 2
        c = bernoulli(2 * s) * w(QQ, 1, 2 * s) / (4 * s)
        return abs(c)


#If E is a totally real Galois extension of Q and Galois group is p-elementary, 
#then given a list of intermediate field of E and a list of associated conductors of these intermediate fields.
#We can compute the size of even K groups of ring of integers of E by the formula in Theorem 4.1.
def MultiSize(list1, list2, n):
    """
    Computes a combined size metric over multiple fields and conditions.
    :param list1: A list of intermediate fields.
    :param list2: A list of conductors corresponding to the intermediate fields.
    :param n: An integer n=2k.
    :return: Size of the K-group for number field E.
    """
    prod = 1
    p = list1[0].degree()  # Degree of the first field.
    k = 2  # Degree from Galois group.

    # Compute the product of sizes for each field in the list.
    for i in range(len(list1)):
        prod *= KSize(list1[i], list2[i], n)

    # Normalize by a scaling factor.
    return prod / (KZ(n)**((p**k - p) // (p - 1)))
\end{lstlisting}\label{pelementary}

\end{document}